\documentclass[preprint,12pt]{elsarticle}

\usepackage{amssymb}
\usepackage{amsmath}
\usepackage{mathtools}
\usepackage{amsthm}

\allowdisplaybreaks[1]
\numberwithin{equation}{section}

\usepackage[a4paper,margin=2.5cm,headheight=15pt,footskip=1cm]{geometry}

\usepackage{microtype}

\biboptions{sort&compress}
\usepackage[colorlinks=true, linkcolor=blue, urlcolor=blue, citecolor=blue]{hyperref}

\newtheorem{theorem}{Theorem}[section]
\newtheorem{proposition}[theorem]{Proposition}
\newtheorem{lemma}[theorem]{Lemma}
\newtheorem{corollary}[theorem]{Corollary}
\newcommand{\R}{\mathbb R}
\newcommand{\ball}{\mathbb B}
\newcommand{\dd}{\,\mathrm d}
\newcommand{\Es}{\mathcal E_s}
\newcommand{\mstar}{m_*(s)}
\newcommand{\mb}{m_{\mathrm b}}
\DeclareMathOperator{\diam}{diam}

\begin{document}

\begin{frontmatter}

\title{The sharp mass threshold for the three-dimensional
	Riesz liquid-drop model near the Coulomb exponent}

\author[inst1]{Zijuan Gao}
\ead{gaozijuan3012@163.com}

\author[inst1]{Chengxiang Zhang}
\ead{zcx@bnu.edu.cn}

\affiliation[inst1]{
	organization={Laboratory of Mathematics and Complex Systems (Ministry of Education), 
		School of Mathematical Sciences, Beijing Normal University},
	city={Beijing},
	postcode={100875},
	country={P. R. China}
}

\begin{abstract}
We determine the exact mass threshold of the three-dimensional
liquid-drop model with repulsive kernel $|x-y|^{-s}$ in a two-sided
neighborhood of the Coulomb exponent $s=1$. The infimum is attained
precisely for $0<m\le\mstar$, with balls as the unique minimizers up to
translation and null sets. Here $\mstar$ is the mass at which one ball and two equal
balls at infinite separation have the same energy. We compute the
first-order variation of this threshold at $s=1$. For $s>1$, a
diameter bound proportional to the mass supplies a lower bound on the
averaged planar slicing integral and yields a uniform mass bound.
The resulting compactness and a uniform quadratic energy
gap for finite-perimeter sets close in measure to a ball allow us to
transfer the Coulomb classification to nearby exponents. At
$m=\mstar$, every minimizing sequence has a subsequence asymptotic
in measure to single balls or pairs of equal balls with diverging
separation, with both alternatives realized by minimizing sequences. For $\mstar<m\le3\mstar/2$,
only the separated pair occurs and the infimum is twice the energy
of a half-mass ball. We also obtain the first-order change in
interaction energy with respect to the exponent along subsequences of near-minimizers as $s\to1$
and $m\to m_*(1)$.
\end{abstract}


\end{frontmatter}

\section{Introduction}\label{sec:intro}

We consider the three-dimensional liquid-drop energy
\begin{equation}\label{eq:energy}
  \Es(E)=P(E)+D_s(E),\qquad
  D_s(E)=\frac12\int_E\int_E|x-y|^{-s}\,\dd y\,\dd x,
  \qquad 0<s<3,
\end{equation}
for measurable sets $E\subset\R^3$ of finite volume. Here
$P(E)$ is the De Giorgi perimeter,
\[
 P(E)=\sup_{\substack{\varphi\in C_c^1(\R^3;\R^3)\\
                         \|\varphi\|_{L^\infty}\le1}}
       \int_E\operatorname{div}\varphi\,\dd x,
\]
and $|E|$ denotes Lebesgue measure. The interaction is finite for
every such set: the kernel is locally integrable for $s<3$ and
bounded outside the unit ball. At fixed volume, balls minimize
perimeter by isoperimetry and maximize $D_s$ by Riesz rearrangement
\cite{Maggi,Riesz,LL}. We write
\[
 e_s(m)=\inf\{\Es(E):|E|=m\}.
\]
The case $s=1$ is Gamow's Coulomb model \cite{Gamow,CMT}.
Let $B_m$ denote the ball of volume $m$ centered at the origin.
Following Frank and Nam \cite{FN}, for each
$s\in(0,3)$ let $\mstar$ denote the unique mass at which
the energy of a ball equals that of two infinitely
separated balls of equal mass:
\[
\Es(B_{\mstar})=2\Es(B_{\mstar/2}).
\]
The standard two-ball comparison \eqref{eq:two-ball-upper}, together with the
explicit scaling of the ball energy, yields
$\mb(s)\le\mstar$, where
$\mb(s):=\sup\{m>0:e_s(m)=\Es(B_m)\}$.
Thus, $\mstar$ is the natural candidate for the
transition mass.

The two-ball comparison raises the question of whether
$m_*(s)$ gives the actual transition between ball
minimality and nonattainment. For general Riesz exponents,
related aspects of this question have been studied in
different mass regimes. At small masses, the main results
concern the existence of minimizers and the minimality
and uniqueness of balls, whereas at larger masses the
focus shifts to nonattainment and the structure of
minimizing configurations. Kn\"upfer and Muratov established existence at
sufficiently small masses and uniqueness of balls for $0<s<2$, while Figalli, Fusco, Maggi,
Millot and Morini proved uniqueness of balls in this regime for all $0<s<3$
\cite{KM,FFMMM}.
Frank and Nam established attainment up to the two-ball
mass and strict binding below it \cite{FN}. Their Theorem~2 also gives ball minimality up to $\mstar$ and uniqueness below it under the assumption of nonattainment for every $m>\mstar$. At sufficiently large masses, nonattainment was proved by Kn\"upfer and Muratov
for $s<2$ \cite{KM} and by Frank and Nam for $s=2$ \cite[Theorem~3]{FN}. Related results include thickness
estimates and nonattainment under strong repulsion
in restricted classes \cite{MPR}, as well as the
construction of minimizing sequences from separated ball
configurations for sufficiently small exponents \cite{BC}.


Bronsard, Merlet and Pegon established connected nonball minimizers for certain truncated Coulomb and Yukawa interactions, including nonradial minimizers for the truncated Coulomb kernel. For the unscreened Riesz kernel, they showed that balls have a lower energy-per-mass ratio than infinite cylinders for integer exponents \(1<s<3\), which in dimension three means \(s=2\). This comparison does not establish the full global classification \cite{BMP}.

For the Coulomb exponent $s=1$, the picture is now sharp.
Uniqueness of balls was established for small masses
by Kn\"upfer and Muratov \cite{KM} and by Julin \cite{Julin} and subsequently extended to $m\le1$ by Chodosh and
Ruohoniemi \cite{CR}. Large-mass nonattainment
was obtained by Lu and Otto, Frank, Killip and Nam, and
Schulz \cite{LO,FKN,Schulz}. Chodosh and Gianocca
completed the classification, proving uniqueness of balls
up to the critical mass $m_c=m_*(1)$ and nonattainment above it
\cite{CG}. A complementary description in terms of
generalized minimizers consisting of separated components
was developed by Kn\"upfer, Muratov and Novaga
\cite{KMN}.

The sharp Coulomb classification raises a natural question:
whether the transition between ball minimizers and
nonattainment persists under perturbations of the Riesz
exponent. More precisely, does the mass determined by
the equal two-ball comparison remain the exact threshold
when the exponent moves away from the Coulomb value?
The Coulomb classification alone does not answer this
question. A perturbation of the exponent affects the
interaction at both short and large distances, requiring
uniform control of the compactness and stability estimates.

We answer this question affirmatively. The two-ball mass
remains the exact threshold for attainment and ball
minimality in a two-sided neighborhood of the Coulomb
exponent. Moreover, uniqueness of the ball persists
up to and including the critical mass.

\begin{theorem}\label{thm:main}
There exists $\varepsilon_0>0$ such that, whenever
$|s-1|<\varepsilon_0$ and $m>0$, the infimum $e_s(m)$ is attained
if and only if $m\le\mstar$. For these masses,
$e_s(m)=\Es(B_m)$, and the minimizers are precisely the
translates of $B_m$, up to null sets. In particular,
$\mb(s)=\mstar$ for $|s-1|<\varepsilon_0$.
\end{theorem}

As a consequence of Theorem~\ref{thm:main} and
Proposition~\ref{prop:slope}, the transition mass has
the expansion
\[
m_{\mathrm b}(s)
=m_c+m_*'(1)(s-1)+O((s-1)^2)
\qquad\text{as }s\to1,
\]
where $m_c\approx3.512072$ and
$m_*'(1)\approx-0.019598$.

We also classify minimizing sequences at the transition
mass and in the range immediately above it. At the
critical mass, both compactness to a single ball and
splitting into two equal balls can occur, whereas above
the critical mass only the latter occurs. In addition,
we obtain a first-order expansion of the interaction
energy for near-minimizers as the exponent approaches
the Coulomb value.

The proof requires uniform control of the nonlocal
interaction as the Riesz exponent approaches the Coulomb
value. The main difficulty arises for $s>1$, where the
weight $|x-y|^{1-s}$ in the planar slicing inequality
decays at large distances and does not directly provide
a uniform lower bound. We overcome this difficulty by
combining the slicing inequality with a diameter estimate,
obtaining uniform bounds on the mass, diameter and
perimeter. These estimates yield compactness up to
translations and force any limiting configuration of
minimizers with masses at or above the two-ball threshold to
be a Coulomb minimizer of mass $m_c$, hence a ball.

A second difficulty is to exclude nonball minimizers
that converge to this limiting ball. Convergence alone
does not rule out such configurations. We resolve this
issue by establishing a second-variation estimate that
is uniform in the exponent near $s=1$ and extending the
resulting local stability to finite-perimeter sets through
a penalization and selection argument. This provides the
uniform local rigidity needed to exclude nonball
minimizers and establish the exact transition threshold.
Scaling then propagates minimality at the critical mass
to all smaller masses. Finally, concentration compactness \cite{Lions}
gives the classification of minimizing sequences, while
uniform estimates for the dependence of the interaction
kernel on $s$ yield the first-order expansion near the
Coulomb exponent.

Section~\ref{sec:scalar} computes the two-ball threshold, proves
a quantitative gap for unequal mass splits, and establishes the
kernel estimates. Section~\ref{sec:local} proves uniform local
stability of the ball. Section~\ref{sec:compact} establishes the
compactness estimates and proves Theorem~\ref{thm:main}.
Section~\ref{sec:profiles} classifies minimizing sequences near
the transition and proves the first-order interaction expansion.

\section{The ball energy and the splitting gap}\label{sec:scalar}

We compute the two-ball threshold $\mstar$ and its derivative at $s=1$.
We then quantify the energy gap for unequal splits at this mass
and establish uniform estimates for the dependence of $D_s$
on the set and on the exponent.

Let $\ball_R(x)$ denote the open ball of radius $R$ centered at $x$.
With $B=\ball_1(0)$ and $\omega=|B|=4\pi/3$, we have $B_m=R_mB$,
where $R_m=(m/\omega)^{1/3}$. For $0<s<3$, the powers of the mass in
$P(B_m)$, $D_s(B_m)$, and $D_s(B_m)/P(B_m)$ are, respectively,
\[
  p=\frac{2}{3},\qquad b_s=\frac{6-s}{3},\qquad q_s=b_s-p=\frac{4-s}{3}.
\]

\begin{proposition}\label{prop:slope}
For $0<s<3$, the infimum satisfies
\begin{equation}\label{eq:two-ball-upper}
	e_s(m)\le2\Es(B_{m/2}),\qquad m>0.
\end{equation}
The unit-ball interaction $d_s\mathrel{=}D_s(B)$ and the
two-ball threshold $\mstar$ satisfy
\begin{equation}\label{eq:threshold}
  d_s=\frac{2^{6-s}\pi^2}{(3-s)(4-s)(6-s)},\qquad
  \mstar=\omega
  \left(\frac{4\pi(2^{1/3}-1)}{d_s(1-2^{(s-3)/3})}\right)^{3/(4-s)} .
\end{equation}
In particular, with $m_c=m_*(1)$,
\begin{align}
  m_c =& 5\,\frac{2^{1/3}-1}{1-2^{-2/3}},\label{eq:mc}\\
  m_*'(1)
  =& m_c\left[
       \frac13\log\frac{m_c}{\omega}+\log 2-\frac{31}{30}
       +\frac{\log 2}{3(2^{2/3}-1)}
     \right].\label{eq:slope}
\end{align}
Consequently $m_*(s)=m_c+m_*'(1)(s-1)+O((s-1)^2)$ as $s\to1$.
\end{proposition}

\begin{proof}
For finite-volume sets $E_1,E_2$ with $|E_1\cap E_2|=0$,
it follows from \eqref{eq:energy} that
\begin{equation}\label{eq:union-interaction}
	D_s(E_1\cup E_2)
	=D_s(E_1)+D_s(E_2)
	+\int_{E_1}\int_{E_2}|x-y|^{-s}\,\dd y\,\dd x.
\end{equation}
If $E_1,E_2\subset\ball_R(0)$ and $|z|>2R$, then $E_1$ and $E_2+z$
are separated by at least $|z|-2R$. Their perimeters add, and
\begin{equation}\label{eq:translated-interaction}
	0\le\int_{E_1}\int_{E_2+z}|x-y|^{-s}\,\dd y\,\dd x
	\le |E_1|\,|E_2|(|z|-2R)^{-s}.
\end{equation}
Testing with $B_{m/2}\cup(B_{m/2}+z)$ and letting $|z|\to\infty$ proves \eqref{eq:two-ball-upper}.

	For $r=|z|\le 2$, the intersection $B\cap(B+z)$ consists of two
	congruent spherical caps of height $1-r/2$, so
	$|B\cap(B+z)|=\pi(16-12r+r^3)/12$.
	Changing variables to $z=x-y$ and using polar coordinates, we obtain
	\begin{align*}
		d_s
		=& \frac12\int_{|z|\le2}|z|^{-s}|B\cap(B+z)|\,\dd z\\
		=& \frac{\pi^2}{6}\int_0^2 r^{2-s}(16-12r+r^3)\,\dd r
		=\frac{2^{6-s}\pi^2}{(3-s)(4-s)(6-s)}.
	\end{align*}
	
	For $t>0$, the dilation $(x,y)\mapsto(tx,ty)$ in \eqref{eq:energy}
	and perimeter scaling give
	\begin{equation}\label{eq:scaling}
	  \begin{gathered}
	    P(tE)=t^2P(E),\qquad D_s(tE)=t^{6-s}D_s(E),\\
	    \Es(B_m)=4\pi R_m^2+d_sR_m^{6-s}\mathrel{=}
	      4\pi(m/\omega)^p+d_s(m/\omega)^{b_s}.
	  \end{gathered}
	\end{equation}
	Evaluating the ball energy at $m$ and $m/2$ in \eqref{eq:scaling}, we find
	\begin{equation}\label{eq:ball-two-difference}
		2\Es(B_{m/2})-\Es(B_m)
		=(m/\omega)^p
		\bigl[4\pi(2^{1/3}-1)
		-d_s(1-2^{1-b_s})(m/\omega)^{q_s}\bigr].
	\end{equation}
	Since $q_s>0$ and $d_s(1-2^{1-b_s})>0$, the bracket decreases
	strictly from $4\pi(2^{1/3}-1)$ to $-\infty$ as $m$ increases
	from $0$ to $\infty$. Its unique zero is the mass $m_*(s)$ in \eqref{eq:threshold}.
	Thus the energy difference has the sign of $m_*(s)-m$.
	
	With $q_1=1$ and $d_1=16\pi^2/15$, \eqref{eq:threshold}
	reduces to \eqref{eq:mc} at $s=1$.
	The equivalent form $m_c=5/(2^{-1/3}+2^{-2/3})$ yields
	$5/2<m_c<5$.
	
	Logarithmic differentiation of \eqref{eq:threshold} gives
	\begin{equation}\label{eq:unit-interaction-log-derivative}
		\frac{d_s'}{d_s}
		=-\log 2+\sum_{j\in\{3,4,6\}}\frac{1}{j-s},
	\end{equation}
	and, using $q_s'=-1/3$,
	\[
	q_s\frac{m_*'(s)}{m_*(s)}
	=\frac13\log\frac{m_*(s)}{\omega}
	-\frac{d_s'}{d_s}
	+\frac{\log 2}{3(2^{b_s-1}-1)}.
	\]
	Substituting $s=1$, $b_1=5/3$ and
	$\left.(d_s'/d_s)\right|_{s=1}=-\log 2+31/30$
	from \eqref{eq:unit-interaction-log-derivative} proves \eqref{eq:slope}.
	
	Since the formula for $m_*$ in \eqref{eq:threshold} is smooth on
	$(0,3)$, Taylor's theorem at $s=1$ proves the expansion.
\end{proof}

We quantify the equal-split comparison of Frank and Nam \cite[Lemma~6]{FN}
by estimating the cost of an unequal split at the critical mass, uniformly for
$0<s<3$.

\begin{lemma}[Energy gap for unequal splits]\label{lem:split}
For $0<s<3$ and $0<t<1$, set
\begin{equation}\label{eq:G}
  \begin{aligned}
    H_s(t)=&1-t^{b_s}-(1-t)^{b_s},\qquad
    A(s)=\frac{2^{1-p}-1}{1-2^{1-b_s}},\\
    G_s(t)=&t^p+(1-t)^p-1-A(s)H_s(t).
  \end{aligned}
\end{equation}
At $m=\mstar$, the splitting cost satisfies
\begin{equation}\label{eq:split-gap-identity}
  \mathcal E_s(B_{tm})+\mathcal E_s(B_{(1-t)m})-\mathcal E_s(B_m)
  =P(B_m)G_s(t).
\end{equation}
Moreover,
\begin{equation}\label{eq:split-gap}
  G_s(t)\ge\frac{17}{144\cdot 2^{18}}\,
                   \min\{t,1-t\}^{2/3}(t-1/2)^2 .
\end{equation}
In particular, $G_s(t)=0$ if and only if $t=1/2$.
\end{lemma}

\begin{proof}
	By \eqref{eq:scaling} and the equal-split balance at $m=\mstar$,
	\[
	  \frac{D_s(B_{\mstar})}{P(B_{\mstar})}
	    =\frac{2^{1-p}-1}{1-2^{1-b_s}}=A(s),\qquad
	  \frac{D_s(B_v)}{P(B_v)}=A(s)(v/\mstar)^{q_s}
	  \quad (v>0).
	\]
	For an arbitrary split of mass $v$, \eqref{eq:G} therefore implies
	\begin{align}
	  &\frac{\Es(B_{tv})+\Es(B_{(1-t)v})-\Es(B_v)}{P(B_v)}
	    \notag\\
	  =&\underbrace{t^p+(1-t)^p-1}_{G_s(t)+A(s)H_s(t)}
	    -\frac{D_s(B_v)}{P(B_v)}
	      \underbrace{\bigl[1-t^{b_s}-(1-t)^{b_s}\bigr]}_{H_s(t)}
	    \label{eq:split-energy}\\
	  =&G_s(t)+A(s)\bigl[1-(v/\mstar)^{q_s}\bigr]H_s(t).
	    \label{eq:gap-any-m}
	\end{align}
	At $v=\mstar$, \eqref{eq:gap-any-m} reduces to \eqref{eq:split-gap-identity}.

	\textit{Step 1. Monotonicity.}
	By symmetry, it suffices to consider $0<t\le1/2$.
	The continuous extension of $G_s$ satisfies $G_s(0)=G_s(1/2)=G_s'(1/2)=0$.
	Differentiating \eqref{eq:G} twice yields
	\begin{align*}
		G_s''(t)
		=A(s)b_s(b_s-1)\bigl[t^{b_s-2}+(1-t)^{b_s-2}\bigr]-p(1-p)\bigl[t^{p-2}+(1-t)^{p-2}\bigr].
	\end{align*}

	For $a=p$ or $a=b_s$,
	\[
	t(1-t)\bigl[t^{a-2}+(1-t)^{a-2}\bigr]
	=\bigl[t^{a-1}+(1-t)^{a-1}\bigr]-\bigl[t^a+(1-t)^a\bigr].
	\]
	Since $0<p<1$ and $1<b_s<2$,
	the derivative on $(0,1/2)$ satisfies
	\[
	  (a-1)\bigl[t^{a-2}-(1-t)^{a-2}\bigr]
	   -a\bigl[t^{a-1}-(1-t)^{a-1}\bigr]
	  \begin{cases}
	    <0,& a=p,\\
	    >0,& a=b_s.
	  \end{cases}
	\]
	Thus $t(1-t)G_s''(t)$ is
	strictly increasing, with leading term $-p(1-p)t^{p-1}$ as $t\to0^+$.

	For $x=b_s-1\in(0,1)$, the bound $2^{v}\le1+v$ on $[0,x]$ implies
    \[
       2^x-1=(\log2)\int_0^x 2^{v}\dd v
       \le (\log2)x(1+x).
    \]
    Using $\log2<3/4$ and $2^{1/3}>5/4$, we obtain
    \begin{align*}
        G_s''(1/2)
        =& 4\left[\frac{(2^{1/3}-1)b_s(b_s-1)}{2^{b_s-1}-1}
        -\frac29\,2^{1/3}\right]\\
        \ge& 4\left[\frac{2^{1/3}-1}{\log2}
             -\frac29\,2^{1/3}\right]
        >4\left[\frac{10}{9}\,2^{1/3}-\frac43\right]>\frac29.
    \end{align*}
    Hence $G_s''$ changes sign once, from negative to positive.
	Since $G_s'(t)\sim p\,t^{p-1}$ as $t\to0^+$ and $G_s'(1/2)=0$,
	$G_s'$ crosses zero once in its decreasing interval.
	Thus $G_s$ first increases and then decreases on $(0,1/2)$.
	
	\textit{Step 2. Uniform lower bound.}
	Since $1-2^{-x}\ge x/2$ by concavity and $1-z^x\le-x\log z$ for $0<z<1$,
    \begin{align*}
       A(s)H_s(t)
       &\le 2(2^{1/3}-1)\bigl[-t\log t-(1-t)\log(1-t)\bigr]\\
       &\le 2(2^{1/3}-1)t\bigl(1+\log(1/t)\bigr).
    \end{align*}
    Here $-(1-t)\log(1-t)\le t$. For $t\le1/2$, $1-(1-t)^p\le p\,2^{1-p}t$, so
    \[
       G_s(t)\ge t^p\left[1-t^{1/3}
       \bigl(p\,2^{1-p}+2(2^{1/3}-1)(1+\log(1/t))\bigr)\right].
    \]
    Since $t^{1/3}(1+\log(1/t))$ increases on $(0,2^{-12}]$,
    for $0<t\le2^{-12}$, we have
    \[
    \frac{G_s(t)}{t^p}\ge 1-\frac1{16}
    \left[\frac89+\frac23(1+12\log 2)\right]>1-\frac{17}{36}>\frac12.
    \]

    Set $h=2^{-9}$. By convexity of $x^{b_s-2}$, the $b_s$ term in $G_s''$ is
	minimal at $t=1/2$. On $[1/2-h,1/2]$, we have $t\ge1/4$ and
	$(1-t)^{p-2}\le2^{2-p}$. The mean value theorem and
	$p(1-p)(2-p)4^{3-p}<10$ imply
	\[
	G_s''(t)\ge G_s''(1/2)-p(1-p)(2-p)4^{3-p}h
	>\frac29-10h>\frac15.
	\]
Integrating twice from $1/2$, we obtain
	\[
	G_s(t)=\int_t^{1/2}(r-t)G_s''(r)\,\dd r
	\ge\frac1{10}(1/2-t)^2.
	\]
	
	By Step~1, the minimum on $[2^{-12},1/2-h]$ is at an endpoint, so
    \[
    G_s(t)\ge \min\{2^{-9},h^2/10\}=h^2/10.
    \]
    Since $h^2=2^{-18}$, $t^p(1/2-t)^2\le1/4$, and
    $4(h^2/10)=(2/5)h^2>(17/144)h^2$, the three
    estimates establish \eqref{eq:split-gap} on $(0,1/2]$.
    By symmetry, the estimate also holds on $(0,1)$.
\end{proof}

For measurable sets $E,F\subset\R^3$, let $\chi_E$ denote the indicator of $E$ and write
$E\triangle F=(E\setminus F)\cup(F\setminus E)$ for their symmetric difference.
Convergence of sets in $L^1$ or $L^1_{\mathrm{loc}}$ means convergence
of their indicators in the corresponding space.

\begin{lemma}[Uniform estimates for the Riesz interaction]\label{lem:kernel-estimates}

For a compact interval $I=[s_-,s_+]\subset(0,3)$ and $M>0$, set
\begin{equation}\label{eq:kernel-moment}
  K_{0,M}=\frac{4\pi}{3-s_+}+M,\qquad
  K_{k,M}=\frac{4\pi k!}{(3-s_+)^{k+1}}
        +M\left(\frac{k}{\mathrm e s_-}\right)^k
  \quad(k=1,2).
\end{equation}
Let $E,F\subset\R^3$ be measurable with $|E|,|F|\le M$. Equation \eqref{eq:kernel-moment} specifies the dependence of $K_{k,M}$ on $I$.
For $s,t\in I$ and $k=0,1,2$,
\begin{equation}\label{eq:kernel-unif}
  \sup_{x\in\R^3}\int_E |x-y|^{-s}
  \bigl|\log|x-y|\bigr|^k\dd y\le K_{k,M}.
\end{equation}
Moreover,
\begin{align}
  |D_s(E)-D_s(F)|
    \le& K_{0,M}|E\triangle F|, \label{eq:riesz-L1}\\
  |(D_s-D_t)(E)-(D_s-D_t)(F)|
    \le& K_{1,M}|s-t|\,|E\triangle F|, \label{eq:parameter-L1}\\
  |D_s(E)-D_t(E)|
    \le& \tfrac12 |E|K_{1,M}|s-t|. \label{eq:exponent-lipschitz}
\end{align}
The map $s\mapsto D_s(E)$ belongs to $C^2((0,3))$, with
\begin{equation}\label{eq:kernel-derivatives}
  \begin{aligned}
    \partial_s^k D_s(E)
    =&\frac12\int_E\int_E
      |x-y|^{-s}(-\log|x-y|)^k\dd y\dd x,
      &&0<s<3,\\
    |\partial_s^k D_s(E)|
    \le& \tfrac12|E|K_{k,M}, &&s\in I,
  \end{aligned}
\end{equation}
for $k=0,1,2$.
\end{lemma}

\begin{proof}[\proofname.]
For $k\ge1$, the function $\mathrm e^{-s_-v}v^k$
attains its maximum at $v=k/s_-$.
Splitting at $|x-y|=1$ and setting $v=|\log r|$ with $r=|x-y|$ in polar coordinates, we obtain
\[
  \begin{aligned}
    &\sup_x\int_E \sup_{\sigma\in I} |x-y|^{-\sigma}
           \bigl|\log|x-y|\bigr|^k\dd y\\
    \le&4\pi\underbrace{\int_0^\infty
         \mathrm e^{-(3-s_+)v}v^k\dd v}_{k!/(3-s_+)^{k+1}}
       +M\sup_{v\ge0}\mathrm e^{-s_-v}v^k=K_{k,M}.
  \end{aligned}
\]
By symmetry,
\[
 D_s(E)-D_s(F)
 =\frac12\iint_{\R^3\times\R^3}
  \frac{(\chi_E-\chi_F)(x)(\chi_E+\chi_F)(y)}{|x-y|^s}\dd x\dd y .
\]
Using $|\chi_E-\chi_F|=\chi_{E\triangle F}$ and \eqref{eq:kernel-unif}
with $k=0$ proves \eqref{eq:riesz-L1}. Subtracting the same identity at $t$
and using $r^{-s}-r^{-t}=-\int_t^s r^{-\widetilde s}\log r\dd\widetilde s$
proves \eqref{eq:parameter-L1} by \eqref{eq:kernel-unif} with $k=1$.

Dominated convergence on each compact subinterval of $(0,3)$
yields the derivative formula in \eqref{eq:kernel-derivatives}. The bound on $I$ follows from \eqref{eq:kernel-unif}.
Integrating the bound for $k=1$ in the exponent yields
\eqref{eq:exponent-lipschitz}.
\end{proof}

\section{Uniform local stability of the ball}\label{sec:local}

Fix $I=[s_-,s_+]$ with $0<s_-<1<s_+<2$.
For a measurable set $F\subset\R^3$ of finite volume, set
$\alpha(F)=\inf_{z\in\R^3}|F\triangle(B+z)|$.
For $0<s<3$ and $m>0$, define
\begin{equation}\label{eq:normalized}
	\mathcal F_{s,m}(F)=P(F)+R_m^{4-s}D_s(F).
\end{equation}
If $|F|=\omega$, then $|R_mF|=m$. By \eqref{eq:scaling},
$\mathcal E_s(R_mF)=R_m^2\mathcal F_{s,m}(F)$.
We first prove a uniform quadratic energy gap for small radial
graphs with the same volume and barycenter as the ball. A
penalization argument extends the estimate to sets of finite
perimeter close to a ball in measure.

\subsection{The spherical quadratic form}

We identify the quadratic part of the energy at the ball and
diagonalize it in spherical harmonics.

Let $\mathcal H^2$ denote two-dimensional Hausdorff measure.
On $\mathbb S^2=\partial B$, let $dS$ denote its restriction and
$\nabla_{\mathbb S^2}$ the tangential gradient. For real-valued
$u,v\in L^2(\mathbb S^2)$, set
$\langle u,v\rangle=\int_{\mathbb S^2}uv\,dS$ and
$\|u\|_2^2=\langle u,u\rangle$. On $H^1(\mathbb S^2)$, use
$\|u\|_{H^1}^2=\|u\|_2^2+\|\nabla_{\mathbb S^2}u\|_2^2$.

For $0<s<2$ and $u\in L^2(\mathbb S^2)$, define
\[
(T_su)(\theta)=\int_{\mathbb S^2}
\frac{u(\xi)}{|\theta-\xi|^s}\,dS_\xi
\qquad\text{for a.e. }\theta\in\mathbb S^2.
\]
We write $\lambda_\ell(s)$ for the eigenvalue of $T_s$ on spherical
harmonics of degree $\ell$.

For $u\in C^1(\mathbb S^2)$ with $1+u>0$, set
\[F_u=\{r\theta:\theta\in\mathbb S^2,\ 0\le r<1+u(\theta)\},\quad 
\|u\|_{C^1}=\|u\|_\infty+\|\nabla_{\mathbb S^2}u\|_\infty.\]
For $0<s<2$ and $m>0$, define the quadratic form on real-valued
$H^1(\mathbb S^2)$ by
\begin{equation}\label{eq:form}
	\begin{aligned}
		Q_{s,m}(u)
		=\int_{\mathbb S^2}
		\bigl(|\nabla_{\mathbb S^2}u|^2-2u^2\bigr)\,dS
		+R_m^{4-s}
		\bigl(\langle u,T_su\rangle-\lambda_1(s)\|u\|_2^2\bigr).
	\end{aligned}
\end{equation}

\begin{lemma}\label{lem:spectrum}
Let $s\in I$, $m>0$
and $0<\eta<1/4$. If $u\in C^1(\mathbb S^2)$ satisfies
$\|u\|_{C^1}\le\eta$ and $|F_u|=\omega$, then
\begin{equation}\label{eq:graph-taylor}
  \left|\mathcal F_{s,m}(F_u)-\mathcal F_{s,m}(B)-\frac12 Q_{s,m}(u)\right|
  \le C_{\mathrm{Tay}}(1+R_m^{4-s})\eta^{2-s_+}\|u\|_{H^1}^2 .
\end{equation}
The constant $C_{\mathrm{Tay}}=C_{\mathrm{Tay}}(I)>0$
is independent of $s,m,\eta,u$.
For $0<s<2$ and each integer $\ell\ge0$,
the eigenvalue of $T_s$ on degree-$\ell$ spherical harmonics is
\begin{equation}\label{eq:eigenvalues}
  \lambda_\ell(s)=\frac{2^{3-s}\pi}{2-s}
                   \frac{(s/2)_\ell}{(2-s/2)_\ell}.
\end{equation}
Here $(a)_\ell=a(a+1)\cdots(a+\ell-1)$ and $(a)_0=1$.
\end{lemma}

\begin{proof}
	\textit{Step 1. Spherical eigenvalues.}
	Taking $\theta$ as the pole and setting $\rho=|\theta-\xi|$, we obtain
	\begin{equation}\label{eq:sphere-kernel}
	  \int_{\mathbb S^2}\frac{dS_\xi}{|\theta-\xi|^s}
	   =2\pi\int_0^2\rho^{1-s}\dd\rho=\frac{2^{3-s}\pi}{2-s}.
	\end{equation}
    By Schur's test and symmetry of the kernel, $T_s$ is bounded and
	self-adjoint on $L^2(\mathbb S^2)$.
	Let $\mathrm P_\ell$ be the degree-$\ell$ Legendre polynomial with
	$\mathrm P_\ell(1)=1$ and put $t=\theta\cdot\xi$. Since $(1-t)^{-s/2}\in L^1(-1,1)$ for $0<s<2$,
	we may apply the Funk--Hecke formula. Rodrigues' formula and
	$\ell$ integrations by parts then yield
	\begin{align*}
		\lambda_\ell(s)
		=& 2^{1-s/2}\pi\int_{-1}^1
		(1-t)^{-s/2}\mathrm P_\ell(t)\,\dd t\\
		=& \frac{2^{1-s/2}\pi(-1)^\ell(s/2)_\ell}{2^\ell\ell!}
		\int_{-1}^1(1-t)^{-s/2-\ell}(t^2-1)^\ell\,\dd t\\
		=& \frac{2^{1-s/2}\pi(s/2)_\ell}{2^\ell\ell!}
		\int_{-1}^1(1-t)^{-s/2}(1+t)^\ell\,\dd t .
	\end{align*}
	
	For $0\le j<\ell$, the boundary terms satisfy
	\[
	\partial_{t}^j(1-t)^{-s/2}\,
	\partial_{t}^{\ell-1-j}(t^2-1)^\ell
	=O\bigl((1-t)^{1-s/2}(1+t)^{j+1}\bigr).
	\]
	Thus they vanish as $t\to\pm1$ since $s<2$.
	
	With $t=2v-1$, the beta integral yields
	\begin{align*}
		\lambda_\ell(s)
		=\frac{2^{2-s}\pi(s/2)_\ell}{\ell!}
		\int_0^1 v^\ell(1-v)^{-s/2}\,\dd v
		=\frac{2^{2-s}\pi(s/2)_\ell}{(1-s/2)_{\ell+1}}.
	\end{align*}
	This proves \eqref{eq:eigenvalues}.
	
	\textit{Step 2. Perimeter expansion.}
	Under the parametrization $x=(1+u(\theta))\theta$, the surface
	element is
	\[
	\dd\mathcal H^2(x)
	=(1+u)\sqrt{(1+u)^2+|\nabla_{\mathbb S^2}u|^2}\,dS_\theta.
	\]
	Thus
	\[
	P(F_u)=\int_{\mathbb S^2}(1+u)^2
	\sqrt{1+\frac{|\nabla_{\mathbb S^2}u|^2}{(1+u)^2}}\,dS .
	\]
	The constraint $\int_{\mathbb S^2}(1+u)^3\,dS=3\omega$ implies
	\[
	\int_{\mathbb S^2}(2u+u^2)\,dS
	=-\int_{\mathbb S^2}\bigl(u^2+\tfrac23u^3\bigr)\,dS .
	\]
	Apply $0\le 1+x/2-\sqrt{1+x}\le x^2/8$ with
	$x=|\nabla_{\mathbb S^2}u|^2/(1+u)^2$. Since $1+u\ge3/4$ and
	$\|u\|_{C^1}\le\eta$,
	\begin{equation}\label{eq:perimeter-taylor}
		\begin{aligned}
			&\left|P(F_u)-P(B)-\frac12\int_{\mathbb S^2}
			(|\nabla_{\mathbb S^2}u|^2-2u^2)\,dS\right|\\
			\le& \frac23\int_{\mathbb S^2}|u|^3\,dS
			+\frac18\int_{\mathbb S^2}
			\frac{|\nabla_{\mathbb S^2}u|^4}{(1+u)^2}\,dS\\
			\le& \frac23\eta\|u\|_2^2
			+\frac29\eta^2\|\nabla_{\mathbb S^2}u\|_2^2 .
		\end{aligned}
	\end{equation}

	\textit{Step 3. Interaction expansion.}
	Let $v_s(|x|)=\int_B|x-y|^{-s}\dd y$ be the radial potential of $B$
	and set $f=\chi_{F_u}-\chi_B$. Expanding \eqref{eq:energy} and using
	$\int_{\R^3}f\dd x=0$, we obtain
	\begin{equation}\label{eq:D-diff}
		\begin{aligned}
			D_s(F_u)-D_s(B)
			=& \underbrace{\int_{\R^3}
				\bigl(v_s(|x|)-v_s(1)\bigr)f(x)\dd x}_{I_1}
			+\underbrace{\frac12\iint_{\R^3\times\R^3}
				\frac{f(x)f(y)}{|x-y|^s}\dd y\dd x}_{I_2}.
		\end{aligned}
	\end{equation}
	
	For $r,t\in[1-\eta,1+\eta]$ and $\theta,\xi\in\mathbb S^2$, write
	$\rho=|\theta-\xi|$ and $\kappa=|r\theta-t\xi|$. Then
	\[ 
	\kappa^2=(r-t)^2+rt\rho^2\ge(3\rho/4)^2,\quad 
	|\kappa-\rho|\le2\eta.
	 \] 
	For $0<\rho\le2$, the mean value theorem and the bounds
	$\kappa^{-s}\le2\rho^{-s}$, $r^2t^2<3$, and $|r^2t^2-1|\le8\eta$
	imply
	\begin{align*}
		|\kappa^{-s}-\rho^{-s}|
		\le& \min\bigl\{\kappa^{-s}+\rho^{-s},\,
		2s(4/3)^{s+1}\eta\rho^{-s-1}\bigr\}\\
		\le& \min\{3\rho^{-s},\,10\eta\rho^{-s-1}\},\\
		|r^2t^2\kappa^{-s}-\rho^{-s}|
		\le& \min\bigl\{(2r^2t^2+1)\rho^{-s},\,
		(10+16\rho)\eta\rho^{-s-1}\bigr\}\\
		\le& \min\{7\rho^{-s},\,42\eta\rho^{-s-1}\}.
	\end{align*}
	Thus
	\[
	\max\bigl\{|\kappa^{-s}-\rho^{-s}|,\,
	|r^2t^2\kappa^{-s}-\rho^{-s}|\bigr\}
	\le g_{s,\eta}(\rho):=
	\begin{cases}
		7\rho^{-s}, & 0<\rho\le\eta,\\
		42\eta\rho^{-s-1}, & \eta<\rho\le2.
	\end{cases}
	\]
	Integrating in spherical coordinates as in \eqref{eq:sphere-kernel}, we obtain
	\begin{equation}\label{eq:row-error}
		\begin{aligned}
			\int_{\mathbb S^2}g_{s,\eta}(|\theta-\xi|)\,dS_\xi
			=& 2\pi\int_0^2g_{s,\eta}(\rho)\rho\dd\rho\\
			\le& 2\pi\left(
			7\int_0^\eta\rho^{1-s_+}\dd\rho
			+42\eta\int_\eta^1\rho^{-s_+}\dd\rho+42\eta
			\right)\\
			=&2\pi\left[\left(\frac7{2-s_+}+\frac{42}{s_+-1}\right)\eta^{2-s_+}
       -\frac{42(2-s_+)}{s_+-1}\eta\right]\\
            \le& C_{\mathrm{row}}\eta^{2-s_+},
			\qquad C_{\mathrm{row}}=2\pi\left(\frac7{2-s_+}+\frac{42}{s_+-1}\right).
		\end{aligned}
	\end{equation}
	
	Since $|x|^{-s}$ has a locally integrable weak gradient in $\R^3$ for $s<2$,
	its convolution with $\chi_B$ is $C^1$. 
	Integrating by parts in $B$ and applying
	the degree-one identity for $T_s$, we find
	\begin{align*}
		v_s'(r)
		=& -\int_B\theta\cdot\nabla_y|r\theta-y|^{-s}\dd y
		=-\int_{\mathbb S^2}\frac{\theta\cdot\xi}{|r\theta-\xi|^s}\,dS_\xi,\\
		v_s'(1)
		=& \lim_{r\to1}v_s'(r)
		=-\int_{\mathbb S^2}\frac{\theta\cdot\xi}{\rho^s}\,dS_\xi
		=-\lambda_1(s).
	\end{align*}

	By \eqref{eq:eigenvalues}, $|\lambda_1(s)|\le
	C_\lambda:=2^{3-s_-}\pi/(2-s_+)$. For $r,t\in[1-\eta,1+\eta]$,
	\eqref{eq:row-error} implies
	\begin{align*}
		|r^2v_s'(t)+\lambda_1(s)|
		\le& r^2|v_s'(t)-v_s'(1)|+|r^2-1|C_\lambda\\
		\le& r^2\int_{\mathbb S^2}g_{s,\eta}(\rho)\,dS_\xi
		+|r^2-1|C_\lambda\\
		\le& 2C_{\mathrm{row}}\eta^{2-s_+}+3C_\lambda\eta.
	\end{align*}
	In polar coordinates,
	\begin{equation}\label{eq:interaction-linear-bound}
		\begin{aligned}
			\left|I_1+\frac{\lambda_1(s)}2\|u\|_2^2\right|
			=& \left|\int_{\mathbb S^2}\int_1^{1+u(\theta)}\int_1^r
			\bigl(r^2v_s'(t)+\lambda_1(s)\bigr)\dd t\dd r\,dS_\theta\right|\\
			\le& \left(C_{\mathrm{row}}+\frac32 C_\lambda\right)
			\eta^{2-s_+}\|u\|_2^2.
		\end{aligned}
	\end{equation}

Since $r^2t^2\kappa^{-s}\le6\rho^{-s}$, Fubini's theorem applies.
By Schur's test and \eqref{eq:row-error},
\begin{equation}\label{eq:interaction-quadratic-bound}
	\begin{aligned}
		&\left|I_2-\tfrac12\langle u,T_su\rangle\right|\\
		=& \frac12\left|
		\iint_{\mathbb S^2\times\mathbb S^2}
		\int_1^{1+u(\theta)}\int_1^{1+u(\xi)}
		\bigl(r^2t^2\kappa^{-s}-\rho^{-s}\bigr)
		\dd t\dd r\,dS_\xi\,dS_\theta\right|\\
		\le& \frac12\iint_{\mathbb S^2\times\mathbb S^2}
		g_{s,\eta}(\rho)|u(\theta)u(\xi)|\,dS_\xi\,dS_\theta\\
		\le& \tfrac12 C_{\mathrm{row}}\eta^{2-s_+}\|u\|_2^2.
	\end{aligned}
\end{equation}

Substituting \eqref{eq:perimeter-taylor},
\eqref{eq:interaction-linear-bound}, and \eqref{eq:interaction-quadratic-bound}
in \eqref{eq:normalized}, we bound the left-hand side of
\eqref{eq:graph-taylor} by
\[
\begin{aligned}
	&\frac23\eta\|u\|_2^2
	+\frac29\eta^2\|\nabla_{\mathbb S^2}u\|_2^2
	+\frac32R_m^{4-s}
	(C_{\mathrm{row}}+C_\lambda)\eta^{2-s_+}\|u\|_2^2\\
	\le& 
	\underbrace{\left(\frac89+\frac32
        (C_{\mathrm{row}}+C_\lambda)\right)}_{C_{\mathrm{Tay}}}(1+R_m^{4-s})
	\eta^{2-s_+}\|u\|_{H^1}^2.
\end{aligned}
\]
\end{proof}

The volume and barycenter constraints control the constant and degree-one modes. A spectral gap on the remaining modes turns \eqref{eq:graph-taylor} into a lower bound for the energy deficit.

\begin{lemma}[Uniform coercivity for radial graphs]\label{lem:graph-coercivity}

For $u\in H^1(\mathbb S^2)$, let $u_\ell$ denote its $L^2$-projection
onto the spherical harmonics of degree $\ell$.
Set $u_{\ge2}=\sum_{\ell\ge2}u_\ell$. For $0<s<3$ and $0\le\gamma<4/7$, define
\begin{equation}\label{eq:sharp-mass-range}
   m_\gamma(s)=\frac{\omega}{8}
   \left[\frac{(4-7\gamma)(4-s)(6-s)}{\pi s}\right]^{3/(4-s)},
   \qquad m_{\mathrm{lin}}(s)=m_0(s).
\end{equation}
For $0<s<2$ and $0<m\le m_{\mathrm{lin}}(s)$,
\begin{equation}\label{eq:sharp-spherical-gap}
   \inf_{u_{\ge2}\ne0}\frac{Q_{s,m}(u_{\ge2})}{\|u_{\ge2}\|_{H^1}^2}
   =\frac47\left[1-\left(\frac{m}{m_{\mathrm{lin}}(s)}\right)^{q_s}\right].
\end{equation}
For $0<\gamma<4/7$, $s\in I$ and $0<m\le m_\gamma(s)$,
\begin{equation}\label{eq:uniform-gap}
  Q_{s,m}(u_{\ge2})\mathrel{\ge}\gamma\|u_{\ge2}\|_{H^1}^2.
\end{equation}
There exists $\eta=\eta(I,\gamma)\in(0,1/4)$ such that, if
$u\in C^1(\mathbb S^2)$ satisfies $|F_u|=\omega$, $\int_{F_u}x\dd x=0$, and $\|u\|_{C^1}\le\eta$,
then
\begin{equation}\label{eq:graph-coercivity}
	\mathcal F_{s,m}(F_u)-\mathcal F_{s,m}(B)
	\ge\frac{\gamma}{4}\|u\|_{H^1}^2
	\ge\frac{\gamma\alpha(F_u)^2}{16\pi(5/4)^4}.
\end{equation}
At $s=1$, as $m$ increases from zero, the first zero among the coefficients of $Q_{1,m}$ with $\ell\ge2$ occurs in degree two at $m=10$.
\end{lemma}

\begin{proof}[\proofname.]
The spherical Laplacian satisfies
$-\Delta_{\mathbb S^2}u_\ell=\ell(\ell+1)u_\ell$.
Hence \eqref{eq:form} reads
\begin{equation}\label{eq:quadratic-eigenvalue}
  \begin{aligned}
    Q_{s,m}(u_\ell)=&\mu_\ell(s,m)\|u_\ell\|_2^2,\\
    \mu_\ell(s,m)=&(\ell-1)(\ell+2)
      -R_m^{4-s}\bigl(\lambda_1(s)-\lambda_\ell(s)\bigr),
      \qquad \ell\ge0.
  \end{aligned}
\end{equation}
By \eqref{eq:eigenvalues},
$\lambda_{j+1}(s)/\lambda_j(s)=(j+s/2)/(j+2-s/2)<1$ for $0<s<2$.
Hence, for $\ell\ge2$,
\begin{equation}\label{eq:gap-ratio}
	\begin{aligned}
		\frac{\lambda_1(s)-\lambda_\ell(s)}{(\ell-1)(\ell+2)}
		=& \frac{2-s}{(\ell-1)(\ell+2)}
		\sum_{j=1}^{\ell-1}\frac{\lambda_j(s)}{j+2-s/2}\\
		\le& \frac{\lambda_1(s)-\lambda_2(s)}{\ell+2}
		\le\frac{\lambda_1(s)-\lambda_2(s)}4 .
	\end{aligned}
\end{equation}
At $s=1$, \eqref{eq:eigenvalues} and \eqref{eq:quadratic-eigenvalue} reduce to
\[
\lambda_\ell(1)=\frac{4\pi}{2\ell+1},\qquad
\mu_2(1,m)=4-\frac{8\pi}{15}\frac m\omega
=4\left(1-\frac m{10}\right).
\]
Since the last bound in \eqref{eq:gap-ratio} is attained only at $\ell=2$,
this proves the assertion at $s=1$.

By \eqref{eq:eigenvalues},
\begin{align*}
	\lambda_1(s)-\lambda_2(s)
	=& \frac{2-s}{3-s/2}\lambda_1(s)
	=\frac{2^{4-s}\pi s}{(4-s)(6-s)},\\
	\frac{\mathrm d}{\mathrm d s}\log\bigl(\lambda_1(s)-\lambda_2(s)\bigr)
	=& \frac4{s(4-s)}+\frac1{6-s}-\log2
	>1-\log2>0 .
\end{align*}
In particular,
\[
   \mu_2(s,m)=4-R_m^{4-s}(\lambda_1(s)-\lambda_2(s))
      =4\left[1-\left(\frac{m}{m_{\mathrm{lin}}(s)}\right)^{q_s}\right].
\]
For $m\le m_{\mathrm{lin}}(s)$ and $\ell\ge2$, %
\eqref{eq:gap-ratio} implies
\[
 \frac{\mu_\ell(s,m)}{1+\ell(\ell+1)}
 \ge\frac{(\ell-1)(\ell+2)}{1+\ell(\ell+1)}\,\frac{\mu_2(s,m)}4
 \ge\frac{\mu_2(s,m)}7.
\]
Equality holds for $\ell=2$. Orthogonality proves \eqref{eq:sharp-spherical-gap},
and \eqref{eq:sharp-mass-range} implies $\mu_2(s,m)\ge7\gamma$ for $m\le m_\gamma(s)$. This proves
\eqref{eq:uniform-gap}. Since $\lambda_1-\lambda_2$ is increasing,
\begin{equation}\label{eq:normalized-coefficient}
  R_m^{4-s}=(m/\omega)^{q_s}
  \le\frac{4-7\gamma}{\lambda_1(s)-\lambda_2(s)}
  \le\frac{4}{\lambda_1(s_-)-\lambda_2(s_-)}=:\Gamma_I,\qquad s\in I,\quad 0<m\le m_\gamma(s).
\end{equation}

Under the hypotheses of \eqref{eq:graph-coercivity}, the volume constraint and
$4\int_{F_u}x\dd x=\int_{\mathbb S^2}\theta(1+u)^4\,dS=0$ imply
\begin{equation}\label{eq:constraints}
	\begin{aligned}
		\int_{\mathbb S^2}u\,dS
		=& -\int_{\mathbb S^2}u^2\left(1+\frac u3\right)\,dS,\\
		\int_{\mathbb S^2}\theta u\,dS
		=& -\int_{\mathbb S^2}\theta u^2
		\left(\frac32+u+\frac{u^2}{4}\right)\,dS .
	\end{aligned}
\end{equation}
The projection onto constants and coordinate functions is
\[
(u_{0}+u_1)(\theta)=\frac1{4\pi}\int_{\mathbb S^2}
(1+3\theta\cdot\xi)u(\xi)\,dS_\xi .
\]
By \eqref{eq:constraints} and $\|u\|_\infty\le\eta<1/4$,
\begin{align*}
	\|u_{0}+u_1\|_{H^1}^2
	=& \|u_0\|_2^2+3\|u_1\|_2^2\\
	=& \frac1{4\pi}\left(
	\left|\int_{\mathbb S^2}u\,dS\right|^2
	+9\left|\int_{\mathbb S^2}\theta u\,dS\right|^2\right)\\
	\le& \frac{(1+\eta/3)^2+9(3/2+\eta+\eta^2/4)^2}{4\pi}\|u\|_2^4
	\le\frac{10}{\pi}\|u\|_2^4
	\le40\eta^2\|u\|_2^2 .
\end{align*}

Since $\lambda_0(s)>\lambda_1(s)$, \eqref{eq:quadratic-eigenvalue} implies $\mu_0(s,m)\ge-2$ and $\mu_1(s,m)=0$. Thus \eqref{eq:uniform-gap} yields
\begin{align*}
	Q_{s,m}(u)
	=& Q_{s,m}(u_{\ge2})
	+\mu_0(s,m)\|u_0\|_2^2\\
	\ge&\gamma\|u_{\ge2}\|_{H^1}^2-2\|u_0\|_2^2\\
	\ge&\gamma\|u\|_{H^1}^2-(2+\gamma)\|u_{0}+u_1\|_{H^1}^2
	\ge\gamma\|u\|_{H^1}^2-40(2+\gamma)\eta^2\|u\|_2^2 .
\end{align*}
It follows from \eqref{eq:graph-taylor} that
\[
\mathcal F_{s,m}(F_u)-\mathcal F_{s,m}(B)
\ge\left(\frac\gamma2-C_{\mathrm{Tay}}(1+\Gamma_I)\eta^{2-s_+}
-20(2+\gamma)\eta^2\right)\|u\|_{H^1}^2 .
\]
In polar coordinates,
\begin{align*}
	\alpha(F_u)
	\le& |F_u\triangle B|
	=\frac13\int_{\mathbb S^2}|(1+u)^3-1|\,dS\\
	\le& (1+\eta)^2\int_{\mathbb S^2}|u|\,dS
	\le\left(\frac54\right)^2\sqrt{4\pi}\,\|u\|_2 .
\end{align*}
Since $\eta^2\le\eta^{2-s_+}\to0$ as $\eta\downarrow0$, choose $\eta\in(0,1/4)$ so that
\[
\bigl[C_{\mathrm{Tay}}(1+\Gamma_I)+20(2+\gamma)\bigr]
\eta^{2-s_+}\le \gamma/4.
\]
The energy bound and the estimate for $\alpha(F_u)$ now imply \eqref{eq:graph-coercivity}.
\end{proof}

\subsection{From radial graphs to finite-perimeter sets}

Starting from a sequence violating the local energy gap, we minimize a penalized functional.
Its minimizers converge to the ball as radial graphs, so \eqref{eq:graph-coercivity} applies.

\begin{proposition}\label{prop:local}
Fix $0<\gamma<4/7$. Let $s\in I$ and $0<m\le m_\gamma(s)$.
There exists $\delta=\delta(I,\gamma)>0$, independent of $s$ and $m$, such
that for every finite-perimeter set $F\subset\R^3$
with $|F|=\omega$ and $\alpha(F)\le\delta$,
\begin{equation}\label{eq:local}
  \mathcal F_{s,m}(F)-\mathcal F_{s,m}(B)\ge c\,\alpha(F)^2,
  \qquad c=\frac{\gamma}{32\pi(5/4)^4}.
\end{equation}
\end{proposition}

\begin{proof}
	\textit{Step 1. Penalization.}
	For $|F|<\infty$, both $\chi_F$ and $\chi_B$ belong to $L^2(\R^3)$,
	so their convolution $z\mapsto |F\cap(B+z)|$ is continuous and
	vanishes at infinity. Hence
	\[
	\alpha(F)=\min_{z\in\R^3}|F\triangle(B+z)|
	=|F|+\omega-2\max_{z\in\R^3}|F\cap(B+z)|.
	\]
	Thus $\alpha(F)=0$ precisely when $F$ coincides with a translate of $B$
	up to a null set. For any measurable set $A$ of finite volume, the
	triangle inequality implies
	\begin{equation}\label{eq:asymmetry-continuity}
	  |\alpha(F)-\alpha(A)|\le |F\triangle A|.
	\end{equation}
	
	Suppose that \eqref{eq:local} fails. For each $j$, there exist
	$s_j\in I$, $0<m_j\le m_\gamma(s_j)$, and a set $F_j$ of finite perimeter
	with $|F_j|=\omega$ such that
	\[
	0<\delta_j:=\alpha(F_j)<1/j,\qquad
	\mathcal F_{s_j,m_j}(F_j)-\mathcal F_{s_j,m_j}(B)<c\delta_j^2.
	\]
	Following the selection principle of \cite{CL} and its use in \cite{BC,FFMMM},
	consider the penalized functional
	\[
	\mathcal G_j(F)=\mathcal F_{s_j,m_j}(F)+\Lambda_0|\alpha(F)-\delta_j|,
	\qquad \Lambda_0:=1+\Gamma_I K_{0,\omega}.
	\]
	
	For $|F|=\omega$, choose $z$ such that $|F\triangle(B+z)|=\alpha(F)$.
	Since $P(F)\ge P(B)$ and $R_{m_j}^{4-s_j}\le\Gamma_I$ by \eqref{eq:normalized-coefficient}, applying
	\eqref{eq:riesz-L1} to $F$ and $B+z$ yields
	\begin{equation}\label{eq:penalized-lower}
	\begin{aligned}
		\mathcal G_j(F)-\mathcal F_{s_j,m_j}(B)
		\ge& \Lambda_0|\alpha(F)-\delta_j|
		-R_{m_j}^{4-s_j}K_{0,\omega}\alpha(F)\\
		\ge& \bigl(\Lambda_0-\Gamma_I K_{0,\omega}\bigr)\alpha(F)
		-\Lambda_0\delta_j
		=\alpha(F)-\Lambda_0\delta_j.
	\end{aligned}
\end{equation}
	Since $\mathcal G_j(F_j)=\mathcal F_{s_j,m_j}(F_j)$ and $c<1$, we may choose
	a minimizing sequence $(F_{j,k})_k$ for $\inf_{|F|=\omega}\mathcal G_j(F)$
	such that $\mathcal G_j(F_{j,k})\le\mathcal F_{s_j,m_j}(B)+\delta_j^2$.
	By translation invariance, we may assume that
	$|F_{j,k}\triangle B|=\alpha(F_{j,k})$.
	Applying \eqref{eq:penalized-lower} to $F_{j,k}$ and
	\eqref{eq:kernel-unif} to $B$, we obtain
	\begin{equation}\label{eq:penalized-bounds}
		\begin{aligned}
			|F_{j,k}\triangle B|
			=& \alpha(F_{j,k})
			\le \Lambda_0\delta_j+\delta_j^2
			\le (\Lambda_0+1)\delta_j,\\
			P(F_{j,k})
			\le& \mathcal G_j(F_{j,k})
			\le 4\pi+R_{m_j}^{4-s_j}
			\underbrace{D_{s_j}(B)}_{\le \omega K_{0,\omega}/2}
			+\delta_j^2
			\le 4\pi+\tfrac12\Gamma_I\omega K_{0,\omega}+1=:P_{\max}.
		\end{aligned}
	\end{equation}

\textit{Step 2. Restoring the volume constraint.}
Choose $\psi\in C^\infty(\R;[0,1])$ equal to $1$ on $(-\infty,0]$
and to $0$ on $[1,\infty)$. Let $\Phi_t$ be the flow of
$X(x)=\psi(|x|-2)x$, and write $L_X=\|DX\|_\infty$. Set
\[
\tau_0=\tfrac13\log(4/3),\qquad
\delta_V=\min\left\{1,\frac{3\omega\tau_0}{2},
\frac{3\omega\,\mathrm e^{-3L_X\tau_0}}
{4\|\operatorname{div}X\|_\infty}\right\}.
\]
By Gronwall's inequality, $\|D\Phi_t\|_\infty\le\mathrm e^{L_X|t|}$.
Since $X(x)=x$ on $\ball_2(0)$, we have $\Phi_t B=\mathrm e^t B$
for $|t|\le\tau_0$.

Let $F$ satisfy $|F\triangle B|\le\delta_V$ and $P(F)\le P_{\max}$,
with $P_{\max}$ as in \eqref{eq:penalized-bounds}.
For bounded $\varphi\in C^1(\R^3)$, changing variables $x=\Phi_t(y)$
and differentiating in $t$ yields
\begin{equation}\label{eq:volume-transport}
	\frac{\mathrm d}{\mathrm dt}\int_{\Phi_t F}\varphi\dd x
	=\int_{\Phi_t F}\operatorname{div}(\varphi X)\dd x.
\end{equation}
For $|t|\le\tau_0$, \eqref{eq:volume-transport} with $\varphi=1$
and the bound $\det D\Phi_t\le\mathrm e^{3L_X|t|}$ imply
\begin{align*}
	\frac{\mathrm d}{\mathrm dt}|\Phi_t F|
	=& 3\omega\mathrm e^{3t}
	+\int_{\R^3}(\chi_{\Phi_t F}-\chi_{\Phi_t B})
	\operatorname{div}X\dd x\\
	\ge& 3\omega\mathrm e^{-3\tau_0}
	-\|\operatorname{div}X\|_\infty
	\mathrm e^{3L_X\tau_0}|F\triangle B|
	\ge \frac{3\omega}{2}.
\end{align*}
Since $\bigl||F|-\omega\bigr|\le\delta_V\le 3\omega\tau_0/2$,
we have $|\Phi_{-\tau_0}F|\le\omega\le|\Phi_{\tau_0}F|$.
Thus there is a unique $\tau\in[-\tau_0,\tau_0]$ such that
\begin{equation}\label{eq:volume-correction-time}
	|\Phi_\tau F|=\omega,\qquad
	|\tau|\le\frac{2}{3\omega}\bigl||F|-\omega\bigr|.
\end{equation}

By the area formula,
$P(\Phi_t F)\le\|D\Phi_t\|_\infty^2P(F)\le\mathrm e^{2L_X|t|}P(F)$.
Integrating \eqref{eq:volume-transport} from $0$ to $t$ and taking
the supremum over $\varphi\in C_c^1(\R^3)$ with $\|\varphi\|_\infty\le1$,
we obtain
\begin{equation}\label{eq:flow-displacement}
	|\Phi_t F\triangle F|
	\le \|X\|_\infty\left|\int_0^t P(\Phi_r F)\dd r\right|
	\le \|X\|_\infty |t|\mathrm e^{2L_X|t|}P(F).
\end{equation}
For the translation flow $x\mapsto x+tz$, the same argument at $t=1$ yields
\begin{equation}\label{eq:translation-displacement}
  |(F+z)\triangle F|\le |z|P(F).
\end{equation}

Since $|F|\le\omega+1$, \eqref{eq:normalized-coefficient},
\eqref{eq:riesz-L1} and \eqref{eq:asymmetry-continuity} imply
\begin{equation}\label{eq:volume-fix}
	\begin{aligned}
		\mathcal G_j(\Phi_\tau F)-\mathcal G_j(F)
		\le& (\mathrm e^{2L_X|\tau|}-1)P(F)
		+\left(\Gamma_I K_{0,\omega+1}+\Lambda_0\right)
		|\Phi_\tau F\triangle F|\\
		\le& C_V\bigl||F|-\omega\bigr|.
	\end{aligned}
\end{equation}
By \eqref{eq:volume-correction-time} and \eqref{eq:flow-displacement},
the last inequality holds uniformly in $j$ and $F$ with
\[
C_V=\frac{2P_{\max}\mathrm e^{2L_X\tau_0}}{3\omega}
\left[2L_X+
\left(\Gamma_I K_{0,\omega+1}+\Lambda_0\right)\|X\|_\infty\right].
\]

\textit{Step 3. Existence for the penalized problem.}
Fix $j$ with $(\Lambda_0+1)\delta_j\le\delta_V$.
By \eqref{eq:penalized-bounds}, the indicators $\chi_{F_{j,k}}$ are bounded in
$BV(\R^3)$, the space of functions of bounded variation.
Local compactness supplies a subsequence along which
$\chi_{F_{j,k}}\to\chi_F$ in
$L^1_{\mathrm{loc}}(\R^3)$ and almost everywhere. Fatou's lemma and
lower semicontinuity of perimeter yield
\[
|F|\le\omega,\qquad |F\triangle B|\le(\Lambda_0+1)\delta_j,
\qquad P(F)\le P_{\max}.
\]
Set $h=\omega-|F|$, so that $0\le h\le |F\triangle B|$.
Since $|F_{j,k}|=\omega$, we obtain by dominated convergence on $F$
\begin{equation}\label{eq:penalized-volume-loss}
	|F_{j,k}\triangle F|=h+2|F\setminus F_{j,k}|\longrightarrow h.
\end{equation}

For an open set $U\subset\R^3$, restrict the test fields in the definition of
$P(E)$ in Section~\ref{sec:intro} to $C_c^1(U;\R^3)$ to obtain the relative perimeter $P(E;U)$.
Suppose $E$ has locally finite perimeter. Then $P(E;U)=|D\chi_E|(U)$,
where $D\chi_E$ is the distributional derivative of $\chi_E$
and $|D\chi_E|$ its total variation measure \cite[Chapter~12]{Maggi}.
For each $x\in\R^3$ and a.e. $r>0$, the Gauss--Green formula for
intersection with a ball \cite[Lemma~15.12, Eq.~(15.14)]{Maggi}
and the coarea formula yield
\begin{equation}\label{eq:ball-slicing}
	\begin{aligned}
		P(E\cap\ball_r(x);U)
		=& P(E;U\cap\ball_r(x))
		+\mathcal H^2(E\cap U\cap\partial\ball_r(x)),\\
		P(E\setminus\ball_r(x);U)
		=& P(E;U\setminus\overline{\ball_r(x)})
		+\mathcal H^2(E\cap U\cap\partial\ball_r(x)),\\
		\frac{\mathrm d}{\mathrm dr}|E\cap U\cap\ball_r(x)|
		=& \mathcal H^2(E\cap U\cap\partial\ball_r(x)).
	\end{aligned}
\end{equation}
The exterior identity follows from
$D\chi_{E\setminus\ball_r(x)}=D\chi_E-D\chi_{E\cap\ball_r(x)}$
and the fact that $|D\chi_E|(\partial\ball_r(x))=0$ for a.e. $r>0$.

At volume $m=|E|>0$, the isoperimetric inequality reads
$P(E)\mathrel{\ge}P(B_m)=C_{\mathrm{iso}}m^p$, where $C_{\mathrm{iso}}=3\omega^{1/3}$.
Suppose $|E|+P(E)<\infty$ and fix $R>0$. Applying this
inequality to $E\setminus\ball_r(0)$ and averaging
\eqref{eq:ball-slicing} over $R<r<2R$ with $U=\R^3$, we obtain
\begin{align*}
	P(E;\ball_R(0))+C_{\mathrm{iso}}|E\setminus\ball_{2R}(0)|^p
	\le& \frac1R\int_R^{2R}
	\bigl(P(E;\ball_r(0))+P(E\setminus\ball_r(0))\bigr)\dd r\\
	=& P(E)+\frac{|E\cap(\ball_{2R}(0)\setminus\ball_R(0))|}{R}
	\le P(E)+\frac{|E|}{R}.
\end{align*}
Take $E=F_{j,k}$. For fixed $R>0$, local convergence implies
\[
|F_{j,k}\setminus\ball_{2R}(0)|
=\omega-|F_{j,k}\cap\ball_{2R}(0)|
\longrightarrow \omega-|F\cap\ball_{2R}(0)|.
\]
By lower semicontinuity of perimeter on $\ball_R(0)$, taking first
the lower limit in $k$ and then letting $R\to\infty$ yields
\begin{equation}\label{eq:missing-perimeter}
	\liminf_{k\to\infty}P(F_{j,k})\ge P(F)+C_{\mathrm{iso}}h^p.
\end{equation}

By Fatou's lemma on $\R^3\times\R^3$,
$D_{s_j}(F)\le\liminf_{k\to\infty}D_{s_j}(F_{j,k})$.
For the penalty, \eqref{eq:asymmetry-continuity} and \eqref{eq:penalized-volume-loss} imply
\[
\liminf_{k\to\infty}|\alpha(F_{j,k})-\delta_j|
\ge |\alpha(F)-\delta_j|-h.
\]
These estimates and \eqref{eq:missing-perimeter} bound the defect from below.
Testing with the volume correction of $F$ in \eqref{eq:volume-fix} bounds it from above. Thus
\begin{equation}\label{eq:penalized-defect}
	C_{\mathrm{iso}}h^p-\Lambda_0h
	\le \inf_{|A|=\omega}\mathcal G_j(A)-\mathcal G_j(F)
	\le C_Vh.
\end{equation}
Since $p=2/3$, any $h>0$ would satisfy
\[
\left(\frac{C_{\mathrm{iso}}}{C_V+\Lambda_0}\right)^3
\le h\le(\Lambda_0+1)\delta_j.
\]
As $\delta_j\to0$, this is impossible for all sufficiently large $j$.
Thus $h=0$ and $F_{j,k}\to F$ in $L^1(\R^3)$ by \eqref{eq:penalized-volume-loss}. By \eqref{eq:penalized-defect},
$A_j:=F$ minimizes $\mathcal G_j$ among sets of volume $\omega$,
with $P(A_j)\le P_{\max}$. Passing to the $L^1$ limit in \eqref{eq:penalized-bounds} and using \eqref{eq:asymmetry-continuity}, we obtain
$|A_j\triangle B|=\alpha(A_j)\le(\Lambda_0+1)\delta_j$.

\textit{Step 4. Perimeter comparison and density bounds.}
Set
\[
\Lambda=\Gamma_I K_{0,\omega+1}+\Lambda_0+C_V,\qquad
r_0=\min\left\{\left(\frac{\delta_V}{2\omega}\right)^{1/3},
\frac1{\Lambda}\right\}.
\]
Fix $j$ large enough that $(\Lambda_0+1)\delta_j\le\delta_V/2$.
For $x\in\R^3$ and $0<r\le r_0$, let $A'$ agree with $A_j$ a.e.\
outside a compact subset of $\ball_r(x)$. By locality of perimeter,
it suffices to consider $P(A')\le P(A_j)\le P_{\max}$. Since
\[
|A'\triangle B|
\le |A'\triangle A_j|+|A_j\triangle B|
\le \omega r^3+(\Lambda_0+1)\delta_j
\le \delta_V,
\]
we may apply \eqref{eq:volume-fix} to $A'$. Since $|\Phi_\tau A'|=\omega=|A_j|$,
the minimality of $A_j$ yields
\[
\mathcal G_j(A_j)\le\mathcal G_j(\Phi_\tau A')
\le\mathcal G_j(A')+C_V\bigl||A'|-\omega\bigr|.
\]
Using \eqref{eq:normalized-coefficient}, \eqref{eq:riesz-L1} with $M=\omega+1$,
and \eqref{eq:asymmetry-continuity}, we obtain
\begin{equation}\label{eq:perimeter-comparison}
	\begin{aligned}
		P(A_j;\ball_r(x))-P(A';\ball_r(x))
		=& P(A_j)-P(A')\\
		\le& \left(\Gamma_I K_{0,\omega+1}+\Lambda_0\right)
		|A'\triangle A_j|
		+C_V\underbrace{\bigl||A'|-\omega\bigr|}_{\le |A'\triangle A_j|}\\
		\le& \Lambda |A'\triangle A_j|.
	\end{aligned}
\end{equation}
Thus $A_j$ is a $(\Lambda,r_0)$-minimizer of perimeter for all
sufficiently large $j$. Writing $E^c=\R^3\setminus E$, we have $P(E^c;U)=P(E;U)$ and
$E^c\triangle F^c=E\triangle F$, so the local comparison also holds
for $A_j^c$.

For either phase $E\in\{A_j,A_j^c\}$, set $v(r)=|E\cap\ball_r(x)|$.
For $0<r<r_0$, compare $E$ with $E\setminus\ball_r(x)$ in the fixed
ball $\ball_{r_0}(x)$. By \eqref{eq:perimeter-comparison},
\eqref{eq:ball-slicing}, and the isoperimetric inequality, for
a.e. $r\in(0,r_0)$,
\begin{equation}\label{eq:phase-cut}
	\begin{aligned}
		P(E;\ball_r(x))-v'(r)
		=& P(E;\ball_{r_0}(x))
		-P(E\setminus\ball_r(x);\ball_{r_0}(x))
		\le \Lambda v(r),\\
		C_{\mathrm{iso}}v(r)^{2/3}
		\le& P(E\cap\ball_r(x))
		=P(E;\ball_r(x))+v'(r)
		\le 2v'(r)+\Lambda v(r).
	\end{aligned}
\end{equation}
Since $v'(r)=\mathcal H^2(E\cap\partial\ball_r(x))\le4\pi r^2$ a.e.\
and $v(r)\le\omega r^3$, it follows from \eqref{eq:phase-cut} that
\begin{equation}\label{eq:local-perimeter-growth}
	P(A_j;\ball_r(x))
	\le 4\pi r^2+\Lambda\omega r^3
	\le (4\pi+\omega)r^2.
\end{equation}
By continuity from below of $|D\chi_{A_j}|$, \eqref{eq:local-perimeter-growth}
holds for every $0<r<r_0$.

Let $x\in\operatorname{spt}(\chi_E\,\mathrm dx)$,
so that $v(r)>0$ for every $r>0$. For $0<r\le r_{0}/2$,
\[
\Lambda v(r)^{1/3}
\le \Lambda\omega^{1/3}r
\le \tfrac12\omega^{1/3}
=\tfrac16 C_{\mathrm{iso}}.
\]
Hence the second line of \eqref{eq:phase-cut} implies, for a.e.\
$r\in(0,r_{0}/2)$,
\begin{equation}\label{eq:density-integration}
	2v'(r)\ge\tfrac12 C_{\mathrm{iso}}v(r)^{2/3}
	\quad\Longrightarrow\quad
	(v^{1/3})'(r)=\frac{v'(r)}{3v(r)^{2/3}}
	\ge\frac{C_{\mathrm{iso}}}{12}.
\end{equation}
Integrate \eqref{eq:density-integration} on $[t,r]$, where $v\ge v(t)>0$.
Since $v(t)\le\omega t^3$, letting $t\downarrow0$ yields
\begin{equation}\label{eq:phase-density}
	|E\cap\ball_r(x)|
	\ge \left(\frac{C_{\mathrm{iso}}r}{12}\right)^3
	=\frac{\omega}{64}\,r^3,
	\qquad 0<r\le r_{0}/2.
\end{equation}

\textit{Step 5. Convergence as radial graphs.}
By \eqref{eq:phase-density} and the Lebesgue density theorem, we may
take $A_j=\R^3\setminus\operatorname{spt}(\chi_{A_j^c}\dd x)$.
For this open representative, the relative isoperimetric inequality
implies $\partial A_j=\operatorname{spt}|D\chi_{A_j}|$.

For $x\in A_j\triangle B$ with $|x|\ne1$, set
$r=\min\{r_{0}/2,\bigl||x|-1\bigr|\}$. Apply \eqref{eq:phase-density}
to $A_j$ for $|x|>1$ and to $A_j^c$ for $|x|<1$. By
\eqref{eq:penalized-bounds}, for all sufficiently large $j$,
\[
\frac{\omega r^3}{64}
\le |(A_j\triangle B)\cap\ball_r(x)|
\le |A_j\triangle B|
\le (\Lambda_0+1)\delta_j
< \frac{\omega r_{0}^3}{512}.
\]
Hence $\bigl||x|-1\bigr|=r\le a_j:=(64|A_j\triangle B|/\omega)^{1/3}$
and $\ball_{1-a_j}(0)\subset A_j\subset\ball_{1+a_j}(0)$.
For the barycenter $\bar x_j$ of $A_j$, this inclusion and \eqref{eq:translation-displacement} imply
\begin{align*}
	|\bar x_j|
	=& \frac1\omega\left|\int_{\R^3}x(\chi_{A_j}-\chi_B)\dd x\right|
	\le \frac{1+a_j}{\omega}|A_j\triangle B|,\\
	|(A_j-\bar x_j)\triangle B|
	\le& |A_j\triangle B|+P(B)|\bar x_j|
	\le (4+3a_j)|A_j\triangle B|.
\end{align*}
We continue to denote the translated minimizer by $A_j$. It has
barycenter zero and satisfies
\begin{equation}\label{eq:selection-shell}
	\ball_{1-a_j-|\bar x_j|}(0)\subset A_j
	\subset\ball_{1+a_j+|\bar x_j|}(0),\qquad
	|A_j\triangle B|\le C\delta_j.
\end{equation}

Here and below in this step, $C$ denotes a constant independent of $j$ and of the boundary point.
Testing against $F_j$ and using \eqref{eq:riesz-L1} and the
isoperimetric inequality, we obtain
\begin{equation}\label{eq:selection-energy}
	\begin{aligned}
		\mathcal G_j(A_j)-\mathcal F_{s_j,m_j}(B)
		\le& \mathcal F_{s_j,m_j}(F_j)-\mathcal F_{s_j,m_j}(B)<c\delta_j^2,\\
		0\le P(A_j)-P(B)
		\le& c\delta_j^2+\Gamma_I K_{0,\omega}|A_j\triangle B|.
	\end{aligned}
\end{equation}
For large $j$, both $\partial A_j$ and $A_j\triangle B$ lie in
$\{1/2<|y|<3/2\}$ by \eqref{eq:selection-shell}. Apply Gauss--Green
to $A_j$ and $B$ with a smooth compactly supported field equal to
$y/|y|$ on this annulus. Since $\operatorname{div}(y/|y|)=2/|y|$,
the measure-theoretic outer unit normal $\nu_j$ on the reduced boundary
$\partial^*A_j$ satisfies
\begin{equation}\label{eq:selection-tilt}
	\begin{aligned}
		\frac12\int_{\partial^*A_j}
		\left|\nu_j(y)-\frac{y}{|y|}\right|^2\dd\mathcal H^2(y)
		=& P(A_j)-P(B)
		+2\int_{\R^3}\frac{\chi_B-\chi_{A_j}}{|y|}\dd y\\
		\le& c\delta_j^2+
		\bigl(\Gamma_I K_{0,\omega}+4\bigr)|A_j\triangle B|
		\le C\delta_j.
	\end{aligned}
\end{equation}

Fix $x\in\partial A_j$ and write $\theta=x/|x|$. Let
$\mathcal C_r(x,\theta)$ be the open cylinder centered at $x$ with
axis $\theta$, radius $r$, and half-height $r$. For
$y\in\partial^*A_j\cap\mathcal C_r(x,\theta)$,
\[
\left|\frac{y}{|y|}-\theta\right|^2
=\frac{|x-y|^2-(|x|-|y|)^2}{|x||y|}
\le 8r^2.
\]
By \eqref{eq:selection-tilt} and \eqref{eq:local-perimeter-growth},
for $0<2r<r_0$,
\begin{align*}
	\operatorname{exc}_j(x,r)
	:=&\frac1{2r^2}\int_{\partial^*A_j\cap\mathcal C_r(x,\theta)}
	|\nu_j-\theta|^2\dd\mathcal H^2\\
	\le& \frac{C\delta_j}{r^2}+8P(A_j;\ball_{\sqrt2r}(x))
	\le \frac{C\delta_j}{r^2}+16(4\pi+\omega)r^2.
\end{align*}
Take $\rho_j=\delta_j^{1/3}$. For large $j$, $4\rho_j<r_0$ and
\begin{equation}\label{eq:selection-excess}
	2\sup_{x\in\partial A_j}\operatorname{exc}_j(x,2\rho_j)
	+2\Lambda\rho_j
	\le C\left(\frac{\delta_j}{\rho_j^2}+\rho_j^2+\rho_j\right)
	\le C\delta_j^{1/3}.
\end{equation}
Since the closure of $\mathcal C_{2\rho_j}(x,\theta)$ lies in $\ball_{4\rho_j}(x)$,
the hypotheses of \cite[Theorem~4.1]{CLM} follow from
\eqref{eq:perimeter-comparison} and \eqref{eq:selection-excess},
with cylindrical excess $2\operatorname{exc}_j(x,2\rho_j)$. Thus
\begin{align*}
	\partial A_j\cap\mathcal C_{\rho_j}(x,\theta)
	=& \{x+z+w_{j,x}(z)\theta:\ z\in\theta^\perp,\ |z|<\rho_j\},\\
	\|\nabla w_{j,x}\|_\infty
	\le& C\bigl(2\operatorname{exc}_j(x,2\rho_j)
	+2\Lambda\rho_j\bigr)^{1/4}
	\le C\delta_j^{1/12},
\end{align*}
where $w_{j,x}\in C^{1,1/4}$ and $w_{j,x}(0)=0$.

By \eqref{eq:phase-density}, the two sides of each graph belong to
opposite phases. Thus $\partial A_j=\partial^*A_j$, and
$\nu_j\cdot\theta$ has a fixed nonzero sign on each graph.
If $\nu_j(x)\cdot\theta<0$, then $|\nu_j-\theta|^2\ge2$ throughout
the graph, so
\[
2\operatorname{exc}_j(x,2\rho_j)
\ge \frac{2\pi\rho_j^2}{(2\rho_j)^2}=\frac{\pi}{2},
\]
contrary to \eqref{eq:selection-excess}. Hence
$\nu_j(x)\cdot x/|x|>0$ on $\partial A_j$.
By \eqref{eq:selection-shell}, each ray starts inside $A_j$ and ends
outside. Its intersections with the compact boundary are transverse,
and at each one the ray exits $A_j$. Thus it crosses exactly once.
The radial projection is therefore a $C^1$ diffeomorphism onto
$\mathbb S^2$.

Write $A_j=F_{u^{(j)}}$. At $x=(1+u^{(j)}(\theta))\theta$,
the two parametrizations have the same tangent plane, so
\[
\nabla_{\mathbb S^2}u^{(j)}(\theta)
=-(1+u^{(j)}(\theta))
\frac{\nu_j(x)-(\nu_j(x)\cdot\theta)\theta}{\nu_j(x)\cdot\theta}
=(1+u^{(j)}(\theta))\nabla w_{j,x}(0).
\]
Consequently, by \eqref{eq:selection-shell},
\[
\|u^{(j)}\|_{C^1}
\le a_j+|\bar x_j|+\tfrac32\sup_{x\in\partial A_j}
\|\nabla w_{j,x}\|_\infty
\le C\bigl(\delta_j^{1/3}+\delta_j^{1/12}\bigr)\longrightarrow0.
\]

\textit{Step 6. Contradiction.}
The graphs from Step~5 have volume $\omega$ and barycenter zero.
Since $\|u^{(j)}\|_{C^1}^{2-s_+}=O(\delta_j^{(2-s_+)/12})\to0$,
\eqref{eq:graph-coercivity} and \eqref{eq:selection-energy} imply,
for all sufficiently large $j$,
\[
2c\,\alpha(A_j)^2+\Lambda_0|\alpha(A_j)-\delta_j|
\le \mathcal G_j(A_j)-\mathcal F_{s_j,m_j}(B)<c\delta_j^2.
\]
Since $\delta_j>0$, the two terms on the left imply
\[
0<1-\frac1{\sqrt2}
\le 1-\frac{\alpha(A_j)}{\delta_j}
\le \frac{c}{\Lambda_0}\delta_j\longrightarrow0,
\]
a contradiction.
\end{proof}

\begin{corollary}[Local stability on compact exponent intervals]\label{cor:sharp-local-ranges}
Let $[a,b]\subset(0,3)$ and $M>0$. With both inequalities understood uniformly in
$a\le s\le b$ and $0<m\le M$, we have
\begin{equation}\label{eq:sharp-rectangle}
 \begin{aligned}
 &\mathcal F_{s,m}(F)-\mathcal F_{s,m}(B)\ge c\,\alpha(F)^2
       &&\Longleftrightarrow\quad M<m_{\mathrm{lin}}(b),\\
 & Q_{s,m}(u_{\ge2})\ge\gamma\|u_{\ge2}\|_{H^1}^2\quad\text{for all }u
       &&\Longleftrightarrow\quad M\le m_\gamma(b)
             \quad(0<\gamma<4/7).
 \end{aligned}
\end{equation}
In the first line, there must exist $c,\delta>0$, independent of $s,m$, such that the estimate holds for every
set $F$ of finite perimeter with $|F|=\omega$ and $\alpha(F)\le\delta$.
In the second line $u\in H^1(\mathbb S^2)$.
For $s\ge2$, the interaction part of $Q_{s,m}$ is understood in the difference
form \eqref{eq:cancelled-interaction-form}.
\end{corollary}

\begin{proof}
For $s<2$, the interaction term $\langle u,T_su\rangle-\lambda_1(s)\|u\|_2^2$ in \eqref{eq:form} equals
\begin{equation}\label{eq:cancelled-interaction-form}
  \frac{2^{4-s}\pi}{4-s}\|u\|_2^2
  -\frac12\iint_{\mathbb S^2\times\mathbb S^2}
      \frac{|u(\theta)-u(\xi)|^2}{|\theta-\xi|^s}\,dS_\xi\,dS_\theta .
\end{equation}
Indeed,
\[
\int_{\mathbb S^2}\frac{1-\theta\cdot\xi}{|\theta-\xi|^s}\,dS_\xi
=\pi\int_0^2\rho^{3-s}\dd\rho=2^{4-s}\pi/(4-s).
\]
Formula \eqref{eq:cancelled-interaction-form} defines the constrained
interaction variation also for $2\le s<3$
\cite[Section~7]{FFMMM}. 
It is finite on $H^1(\mathbb S^2)$. Indeed, for $0<\rho=|\theta-\xi|\le2$ one has $\rho^{-s}\le2^{3-s}\rho^{-3}$. Consequently, the fractional Sobolev seminorm and the
	embedding $H^1(\mathbb S^2)\hookrightarrow H^{1/2}(\mathbb S^2)$ give
	\[
	\begin{aligned}
		\iint_{\mathbb S^2\times\mathbb S^2}
		\frac{|u(\theta)-u(\xi)|^2}{|\theta-\xi|^s}
		\,dS_\xi dS_\theta
		\le C\|u\|_{H^{1/2}(\mathbb S^2)}^2
		\le C\|u\|_{H^1(\mathbb S^2)}^2.
	\end{aligned}
	\]

In the notation of \cite[Theorem~7.1 and Proposition~7.4]{FFMMM},
$V_{3-s}=2D_s$ and $\beta=R_m^{4-s}/2$. By \eqref{eq:threshold}, their critical coefficient for classical perimeter on the unit ball becomes
\begin{equation}\label{eq:critical-interaction-coefficient}
 \beta_*(s)=
 \frac{4P(B)}{2s(3-s)d_s}
 =\frac{(4-s)(6-s)}{2^{3-s}\pi s},
 \qquad
 \frac{\beta}{\beta_*(s)}=\left(\frac{m}{m_{\mathrm{lin}}(s)}\right)^{q_s}.
\end{equation}
To compute the degree-two coefficient directly, restrict both angular integrals to $\rho>\varepsilon$ and apply the Funk--Hecke formula. With $t=1-\rho^2/2$ and $\mathrm P_2(t)=(3t^2-1)/2$, the limit as $\varepsilon\downarrow0$ gives
\[
 2\pi\int_{-1}^{1}\frac{t-\mathrm P_2(t)}{(2-2t)^{s/2}}\dd t
  =2\pi\int_0^2\!\left(\rho^{3-s}-\tfrac38\rho^{5-s}\right)\dd\rho
  =\frac{2^{4-s}\pi s}{(4-s)(6-s)}.
\]
Consequently, for every $0<s<3$ and $m>0$,
\[
 Q_{s,m}(u_2)=
 \left(4-R_m^{4-s}\frac{2^{4-s}\pi s}{(4-s)(6-s)}\right)\|u_2\|_2^2
 =4(1-\beta/\beta_*)\|u_2\|_2^2.
\]
At $\beta=\beta_*$ the constrained form is nonnegative and vanishes on degree
two. By affine dependence on $\beta$, for $0<\beta\le\beta_*$ and $\ell\ge2$,
\[
 Q_{s,m}(u_\ell)\ge (1-\beta/\beta_*)[\ell(\ell+1)-2]\|u_\ell\|_2^2
 \ge \tfrac47(1-\beta/\beta_*)\|u_\ell\|_{H^1}^2.
\]
Both bounds are equalities at $\ell=2$. This proves \eqref{eq:sharp-spherical-gap}
for $0<s<3$.

For $0<m\le m_\gamma(s)$, \eqref{eq:critical-interaction-coefficient}
implies $1-\beta/\beta_*\ge7\gamma/4$. By \cite[(8.2) and Remark~8.2]{FFMMM},
volume- and barycenter-constrained graphs satisfy
\[
 \mathcal F_{s,m}(F_u)-\mathcal F_{s,m}(B)
    \ge c_0\frac{7\gamma}{4}\|u\|_{H^1}^2
 \quad\text{if}\quad \|u\|_{C^1}<
 \min\{1/4,(7\gamma/8)\eta_0\},
\]
where $c_0,\eta_0>0$ depend only on the dimension.
By \eqref{eq:critical-interaction-coefficient}, for $s\in[a,b]$, we have
	\[ R_m^{4-s}\le2\max_{a\le t\le b}\beta_*(t)=:\Gamma_{[a,b]}<\infty .\]
This bound and \eqref{eq:riesz-L1} on $[a,b]$ give uniform constants in Steps~1--5
of the proof of Proposition~\ref{prop:local}.
In those steps $\Gamma_{[a,b]}$ bounds the coefficient of $D_s$, and the graphs from Step~5 satisfy the preceding $C^1$ smallness condition for all large $j$.
Since $\alpha(F_u)^2\le4\pi(5/4)^4\|u\|_2^2$, Step~6 applies with
$c=\min\{1/2,7c_0\gamma/[32\pi(5/4)^4]\}$ for $0<m\le m_\gamma(s)$.

For the converse, take $Y(\theta)=\theta_3^2-1/3$ and set
\[
  1+u_\varepsilon(\theta)
   =\frac{1+\varepsilon Y(\theta)}
    {\bigl((4\pi)^{-1}\int_{\mathbb S^2}(1+\varepsilon Y)^3\,dS\bigr)^{1/3}} .
\]
Then $|F_{u_\varepsilon}|=\omega$ and $\int_{F_{u_\varepsilon}}x\dd x=0$ by evenness and the choice of scale.
Since $Y$ is a degree-two harmonic and $u_\varepsilon=\varepsilon Y+O_{C^1}(\varepsilon^2)$,
the second variation formula \cite[Sections~6--7]{FFMMM} yields
\[
 \mathcal F_{s,m}(F_{u_\varepsilon})-\mathcal F_{s,m}(B)
  =2\varepsilon^2(1-\beta/\beta_*)\|Y\|_2^2+o(\varepsilon^2).
\]
The spherical-cap formula in the proof of Proposition~\ref{prop:slope}
shows that $|B\triangle(B+z)|$ is nondecreasing in $|z|$. By central symmetry,
\[
 |B\triangle(B+z)|\le |(B-z)\triangle(B+z)|
      \le 2|F_{u_\varepsilon}\triangle(B+z)|.
\]
Thus
\[
  3\alpha(F_{u_\varepsilon})\ge |F_{u_\varepsilon}\triangle B|
    =|\varepsilon|\int_{\mathbb S^2}|Y|\,dS+O(\varepsilon^2).
\]
Here $\int_{\mathbb S^2}|Y|\,dS=4\pi\int_0^1|t^2-1/3|\dd t=16\pi/(9\sqrt3)$. Thus, for small $\varepsilon\ne0$, $\alpha(F_{u_\varepsilon})$ is bounded above and below by positive multiples of $|\varepsilon|$.
At $m=m_{\mathrm{lin}}(s)$ the deficit is $o(\alpha(F_{u_\varepsilon})^2)$ as $\varepsilon\to0$.

Using $\log v\le v-1$ at $v=(4-s)/s$, we find
\begin{align*}
 \frac{(4-s)^2}{3}\frac{m_\gamma'(s)}{m_\gamma(s)}
 =&\log\frac{(4-7\gamma)(4-s)(6-s)}{\pi s}
   -1-\frac{4-s}{6-s}-\frac{4-s}{s}\\
 \le&\log\frac{(4-7\gamma)(6-s)}{\pi}-2-\frac{4-s}{6-s}
 <\log 8-\frac73<0 .
\end{align*}
Since $m_\gamma(b)\to m_{\mathrm{lin}}(b)$ as $\gamma\downarrow0$, $M<m_{\mathrm{lin}}(b)$ implies $M\le m_\gamma(b)\le m_\gamma(s)$ for some $\gamma>0$ and all $s\in[a,b]$.
Conversely, if $M\ge m_{\mathrm{lin}}(b)$, the degree-two perturbation at $s=b$,
$m=m_{\mathrm{lin}}(b)$ rules out a uniform quadratic deficit.
At $(s,m)=(b,M)$, the spectral inequality with $u=u_2\ne0$
takes the form
\[
 4\left[1-\left(M/m_{\mathrm{lin}}(b)\right)^{q_b}\right]\ge7\gamma
 \quad\Longleftrightarrow\quad M\le m_\gamma(b).
\]
The reverse implication in the second equivalence follows from \eqref{eq:sharp-spherical-gap} and the monotonicity of $m_\gamma$.

\end{proof}

\section{Compactness and the global threshold}\label{sec:compact}

We bound the mass and diameter of every minimizer with
$m\ge m_{\min}>0$, uniformly for $s\in I$. These bounds allow passage to the Coulomb limit.
Local stability at the limiting ball then rules out nonball minimizers at and above the two-ball threshold.

\begin{proposition}\label{prop:mass-diameter}
Fix $m_{\min}>0$. Let $s\in I$, $m\ge m_{\min}$, and let $E$
attain $e_s(m)$. There exist $M,L_{\mathrm{diam}},C_0>0$
depending only on $I$ and $m_{\min}$ such that
\[
 m\le M,\qquad
 \diam\bigl(\operatorname{spt}(\chi_E\,\mathrm dx)\bigr)\le L_{\mathrm{diam}},\qquad
 P(E)\le C_0m.
\]
\end{proposition}

\begin{proof}
\textit{Step 1. Energy and density estimates.}
We adapt the density estimate of Kn\"upfer and Muratov \cite[Lemmas~4.3 and~7.2]{KM}
to obtain a radius independent of $s\in I$.

Let $k=\lceil m\rceil$, so that $1\le k/m\le1+1/m_{\min}$. As in the proof of
\cite[Theorem~3.4]{KM}, compare with $k$ balls of volume $m/k$
whose pairwise distances tend to infinity. By
\eqref{eq:translated-interaction} and \eqref{eq:kernel-unif},
\begin{equation}\label{eq:linear-energy}
	\begin{aligned}
		\frac{P(E)}{m}\le\frac{e_s(m)}{m}
		\le& \frac{k}{m}\Es(B_{m/k})
		\le C_{\mathrm{iso}}\left(\frac{k}{m}\right)^{1/3}
		+\frac{K_{0,1}}2\\
		\le& (1+1/m_{\min})^{1/3}C_{\mathrm{iso}}+\frac{K_{0,1}}2\\
        =&3\bigl[\omega(1+1/m_{\min})\bigr]^{1/3}
          +\frac{2\pi}{3-s_+}+\frac12=:C_0.
	\end{aligned}
\end{equation}

Fix $x\in\R^3$ and set $w=w(r)=|E\cap\ball_r(x)|$. If $w\le m/2$, dilation by $(1-w/m)^{-1/3}$ restores the volume of
$E\setminus\ball_r(x)$ to $m$. Since $q_s>0$ and $D_s(E\setminus\ball_r(x))\le D_s(E)$,
minimality, \eqref{eq:scaling}, and \eqref{eq:ball-slicing} imply,
for a.e.\ such $r$,
\begin{align*}
	(1-w/m)^{b_s}e_s(m)
	\le& (1-w/m)^{q_s}P(E\setminus\ball_r(x))+D_s(E\setminus\ball_r(x))\\
	\le& e_s(m)-P(E;\ball_r(x))+w'(r).
\end{align*}
Since $1<b_s\le b_{s_-}$, Bernoulli's inequality and
\eqref{eq:linear-energy} yield
\begin{equation}\label{eq:cut-density}
	\begin{aligned}
		P(E;\ball_r(x))
		\le& w'(r)+\bigl[1-(1-w/m)^{b_s}\bigr]e_s(m)\\
		\le& w'(r)+b_s\frac{e_s(m)}m\,w
		\le w'(r)+b_{s_-} C_0w.
	\end{aligned}
\end{equation}

Take
\[
 r_d=\min\left\{\left(\frac{m_{\min}}{2\omega}\right)^{1/3},
              \frac{3}{2b_{s_-}C_0}\right\}.
\]
Then $w(r)\le\omega r^3\le m_{\min}/2\le m/2$ for $0<r\le r_d$.
By \eqref{eq:ball-slicing}, \eqref{eq:cut-density}, and the
isoperimetric inequality, for a.e.\ $r\in(0,r_d)$,
\begin{align*}
	2w'(r)
	\ge& P(E\cap\ball_r(x))-b_{s_-} C_0w\\
	\ge& \bigl(C_{\mathrm{iso}}-b_{s_-} C_0w^{1/3}\bigr)w^{2/3}\\
	\ge& \bigl(C_{\mathrm{iso}}-b_{s_-} C_0\omega^{1/3}r_d\bigr)w^{2/3}
	\mathrel{\ge}\tfrac12 C_{\mathrm{iso}}w^{2/3}.
\end{align*}
If $x\in\operatorname{spt}(\chi_E\dd x)$, then $w(t)>0$ for $t>0$.
Integrating as in \eqref{eq:density-integration} on $[t,r]$ and
letting $t\downarrow0$, we obtain
\begin{equation}\label{eq:density}
	|E\cap\ball_r(x)|
	\ge \left(\frac{C_{\mathrm{iso}}r}{12}\right)^3
	=\frac{\omega}{64}r^3,\qquad 0<r\le r_d.
\end{equation}

\textit{Step 2. Diameter bound.}
Let $\{\ball_{r_d}(x_i)\}$ be a maximal disjoint family with centers in
$\operatorname{spt}(\chi_E\dd x)$. By \eqref{eq:density} and maximality, we have
\begin{equation}\label{eq:diam-packing}
	\begin{gathered}
		N\frac{\omega}{64}r_d^3\le \sum_i |E\cap\ball_{r_d}(x_i)|\le m,
		\quad N:=\#\{x_i\},\quad
		\operatorname{spt}(\chi_E\dd x)
		\subset \bigcup_{i=1}^N\overline{\ball_{2r_d}(x_i)}.
	\end{gathered}
\end{equation}

If this union were disconnected, its components would split $E$,
up to null sets, into two positively separated parts $E_1,E_2$ with
$|E_i|\ge \frac{\omega}{64}r_d^3$. Comparing $E$ with $E_1\cup(E_2+z)$ and using
\eqref{eq:union-interaction} and \eqref{eq:translated-interaction},
\begin{align*}
	0
	\le& \lim_{|z|\to\infty}
	\bigl[\Es(E_1\cup(E_2+z))-\Es(E)\bigr]\\
	=& \underbrace{P(E_1)+P(E_2)-P(E)}_{=\,0\text{ by locality}}
	-\int_{E_1}\int_{E_2}|x-y|^{-s}\dd y\dd x<0.
\end{align*}

Thus any two covering balls are joined by a chain of at most $N$
balls, with consecutive centers at distance at most $4r_d$.
By \eqref{eq:diam-packing},
\begin{equation}\label{eq:diam-linear}
	\begin{aligned}
		\diam\bigl(\operatorname{spt}(\chi_E\dd x)\bigr)
		\le& 2r_d+4r_d(N-1)+2r_d=4r_dN\\
		\le& C_{\mathrm{diam}}m,\qquad
		C_{\mathrm{diam}}:=\max\{1,\frac{256}{\omega r_d^2}\}.
	\end{aligned}
\end{equation}
\textit{Step 3. Planar slicing and the mass bound.}
Fix $\nu\in\mathbb S^2$. Separating the portions of $E$ on opposite sides of
$\{x\cdot\nu=t\}$, as in \cite[Lemma~7]{FN}, costs twice the area of the slice.
Minimality therefore implies, for a.e.\ $t\in\R$,
\begin{equation}\label{eq:plane-slice-cost}
	\begin{aligned}
		\iint_{E\times E}
		\frac{\chi_{\{y\cdot\nu<t<x\cdot\nu\}}}{|x-y|^s}\dd y\dd x
		\le& \sum_{\sigma=\pm1}
		P\bigl(E\cap\{\sigma(x\cdot\nu-t)>0\}\bigr)-P(E)\\
		=& 2\mathcal H^2(E\cap\{x\cdot\nu=t\}).
	\end{aligned}
\end{equation}
For $x,y\in\R^3$,
\begin{align*}
	\int_{\mathbb S^2}\int_\R
	\chi_{\{y\cdot\nu<t<x\cdot\nu\}}\dd t\,dS_\nu
	=& \int_{\mathbb S^2}\max\{(x-y)\cdot\nu,0\}\,dS_\nu\\
	=& 2\pi|x-y|\int_0^1 u\dd u=\pi|x-y|.
\end{align*}
We now integrate \eqref{eq:plane-slice-cost} with respect to $t$ and $\nu$ by Tonelli's theorem and obtain
\begin{equation}\label{eq:planar}
	\iint_{E\times E}|x-y|^{1-s}\dd y\dd x
	\le \frac2\pi\int_{\mathbb S^2}
	\underbrace{\int_\R\mathcal H^2(E\cap\{x\cdot\nu=t\})\dd t}
	_{=\,m\text{ by coarea}}
	\,dS_\nu
	=8m.
\end{equation}

For $s<1$, restrict \eqref{eq:planar} to $|x-y|\ge1$. For $s>1$,
use $|x-y|\le C_{\mathrm{diam}}m$ from \eqref{eq:diam-linear}. Thus
\begin{align*}
	m(m-\omega)
	\le& \int_E|E\setminus\ball_1(x)|\dd x\le8m
	&&(s<1),\\
	m^2\le& 8m &&(s=1),\\
	m^2(C_{\mathrm{diam}}m)^{1-s}
	\le& 8m &&(s>1).
\end{align*}
Hence $m\le8+\omega$ for $s\le1$. For $1<s\le s_+$ and $m\ge1$,
\[
m^{2-s_+}\le m^{2-s}\le8C_{\mathrm{diam}}^{s-1}
\le8C_{\mathrm{diam}}^{s_+-1}.
\]
For $m<1$ the mass bound is immediate. The perimeter and diameter estimates follow from
\eqref{eq:linear-energy} and \eqref{eq:diam-linear} with
\[
M=\max\bigl\{8+\omega,(8C_{\mathrm{diam}}^{s_+-1})^{1/(2-s_+)}\bigr\},
\qquad L_{\mathrm{diam}}=C_{\mathrm{diam}}M.
\]
\end{proof}

\begin{proof}[Proof of Theorem~\ref{thm:main}]
\textit{Step 1. Passage to the Coulomb problem.}
Suppose there exist $s_j\to1$ and nonball minimizers $E_j$ of
$\mathcal E_{s_j}$ at masses $m_j\ge m_*(s_j)$.
Since $m_*(s_j)\to m_c>1$ by Proposition~\ref{prop:slope},
we may assume that $s_j\in I$ and $m_j\ge1$.
By Proposition~\ref{prop:mass-diameter} with $m_{\min}=1$, $m_j\le M$,
$P(E_j)\le C_0M$, and, after translation,
$E_j\subset\overline{\ball_{L_{\mathrm{diam}}}(0)}$ up to null sets.
By $BV$ compactness, along a subsequence,
\[
\chi_{E_j}\to\chi_E \quad\text{in }L^1(\R^3),\qquad
m_j\to m:=|E|\in[m_c,M].
\]
By \eqref{eq:riesz-L1} and \eqref{eq:exponent-lipschitz},
\begin{equation}\label{eq:minimizer-interaction-limit}
	|D_{s_j}(E_j)-D_1(E)|
	\le K_{0,M}|E_j\triangle E|
	+\frac{m}{2}K_{1,M}|s_j-1|
	\longrightarrow 0.
\end{equation}
Lower semicontinuity of perimeter therefore implies
$\mathcal E_1(E)\le\liminf_{j\to\infty}e_{s_j}(m_j)$.

Let $F$ be any set of finite perimeter with $|F|=m$.
The dilate $(m_j/m)^{1/3}F$ has volume $m_j$, so by
\eqref{eq:scaling} and \eqref{eq:exponent-lipschitz},
\begin{align*}
	e_{s_j}(m_j)
	\le& \mathcal E_{s_j}\bigl((m_j/m)^{1/3}F\bigr)\\
	=& \left(\frac{m_j}{m}\right)^{2/3}P(F)
	+\left(\frac{m_j}{m}\right)^{b_{s_j}}D_{s_j}(F)
	\longrightarrow \mathcal E_1(F).
\end{align*}
Taking the infimum over $F$ shows that
\[
 e_1(m)\le\mathcal E_1(E)\le\liminf_j e_{s_j}(m_j)
 \le\limsup_j e_{s_j}(m_j)\le e_1(m).
\]

\textit{Step 2. Rigidity at and above the threshold.}
Chodosh and Gianocca \cite[Theorem~1.1, (1.2), and (1.4)]{CG}
use the normalization \eqref{eq:energy} at $s=1$. By \eqref{eq:mc},
their critical mass is
\[
5\,\frac{2-2^{2/3}}{2^{2/3}-1}
=5\,\frac{2^{1/3}-1}{1-2^{-2/3}}=m_c.
\]
Since $E$ minimizes $\mathcal E_1$ at mass $m\ge m_c$, their
classification implies that $m=m_c$ and $E$ is a ball. After a
common translation of $E_j$ and $E$, we have $E=B_{m_c}$ and
\[
\alpha(R_{m_j}^{-1}E_j)
\le \frac{\omega}{m_j}|E_j\triangle B_{m_j}|
\le \frac{\omega}{m_j}
\bigl(|E_j\triangle B_{m_c}|+|m_j-m_c|\bigr)
\longrightarrow 0.
\]
By \eqref{eq:mc}, \eqref{eq:sharp-mass-range} and the monotonicity of $m_{\mathrm{lin}}$ proved in Corollary~\ref{cor:sharp-local-ranges},
\[
 m_j\to m_c<5<\frac{32}{3\sqrt{\pi}}=m_{\mathrm{lin}}(2)<m_{\mathrm{lin}}(s_+).
\]
Since $|R_{m_j}^{-1}E_j|=\omega$, Corollary~\ref{cor:sharp-local-ranges}
with $[a,b]=I$ and $M=5$ applies for all large $j$.
By \eqref{eq:normalized} and the minimality of $E_j$,
\[
c\,\alpha(R_{m_j}^{-1}E_j)^2
\le R_{m_j}^{-2}
\bigl[e_{s_j}(m_j)-\mathcal E_{s_j}(B_{m_j})\bigr]\le0.
\]
Thus $E_j$ is a ball, contrary to its choice. Hence there exists
$0<\varepsilon_0<\min\{1-s_-,s_+-1\}$ such that, whenever $|s-1|<\varepsilon_0$,
every minimizer of mass $m\ge\mstar$ is a ball.
Fix such an $s$.
Frank and Nam \cite[Theorem~1]{FN} prove attainment at $m=\mstar$,
so $e_s(\mstar)=\Es(B_{\mstar})$.

\textit{Step 3. Smaller masses and nonattainment above the threshold.}
Let $E$ have finite perimeter and volume $0<m<\mstar$, and set
$t=m/\mstar$. Then $|t^{-1/3}E|=\mstar$ and
$t^{-1/3}B_m=B_{\mstar}$. Applying \eqref{eq:scaling}, we obtain
\begin{equation}\label{eq:propagate}
	\begin{aligned}
		\Es(E)-\Es(B_m)
		=&t^{b_s}\underbrace{\bigl[\Es(t^{-1/3}E)-e_s(\mstar)\bigr]}_{\ge0}
		+(1-t^{q_s})\bigl[P(E)-P(B_m)\bigr]\\
		\ge& (1-t^{q_s})\bigl[P(E)-P(B_m)\bigr]\ge0.
	\end{aligned}
\end{equation}
Since $1-t^{q_s}>0$, any minimizer satisfies $P(E)=P(B_m)$
and hence agrees almost everywhere with a translate of $B_m$.

For $m>\mstar$, \eqref{eq:two-ball-upper} and
\eqref{eq:ball-two-difference} imply
\[
e_s(m)\le 2\Es(B_{m/2})<\Es(B_m).
\]
By Step~2, the infimum is not attained. It follows that
$\{m>0:e_s(m)=\Es(B_m)\}=(0,\mstar]$, so $\mb(s)=\mstar$.
\end{proof}

\section{Minimizing sequences at the transition}\label{sec:profiles}

We classify minimizing sequences near the two-ball threshold by comparing the masses of their profiles.
For varying exponents, the limiting profiles determine the first-order change in the interaction.

\subsection{Extraction of profiles}

Countable subadditivity allows us to compare the original mass with any family of profiles extracted by local compactness. Frank and Lieb \cite{FL} used the separation of missing mass to prove compactness under strict binding.

\begin{lemma}[Countable subadditivity]\label{lem:countable-subadd}

Fix $0<s<3$ and set $e_s(0)=0$. For $m_i\ge0$ with $\sum_i m_i<\infty$,
\begin{equation}\label{eq:countable-subadd}
  e_s\!\left(\sum_i m_i\right)\le\sum_i e_s(m_i).
\end{equation}
\end{lemma}

\begin{proof}[\proofname.]
Let $E$ have finite perimeter and volume $m\in(0,\infty)$. By coarea
and \eqref{eq:ball-slicing}, we may choose $r_\ell\in(\ell,\ell+1)$
such that $E_\ell=E\cap\ball_{r_\ell}(0)$ satisfies
\[
0\le P(E_\ell)-P(E;\ball_{r_\ell}(0))
=\mathcal H^2(E\cap\partial\ball_{r_\ell}(0))
\le |E\setminus\ball_\ell(0)|\longrightarrow0.
\]
Thus $P(E_\ell)\to P(E)$. Since $E_\ell\uparrow E$ and $m>0$,
$|E_\ell|>0$ for all sufficiently large $\ell$.
Monotone convergence and \eqref{eq:scaling} imply
\[
\bigl|(m/|E_\ell|)^{1/3}E_\ell\bigr|=m,\qquad
\Es\bigl((m/|E_\ell|)^{1/3}E_\ell\bigr)\to\Es(E).
\]
Hence the infimum defining $e_s(m)$ may be taken over bounded sets.

For bounded competitors $F$ and $G$, \eqref{eq:union-interaction}
and \eqref{eq:translated-interaction} imply
\[
e_s(|F|+|G|)
\le \lim_{|z|\to\infty}\Es\bigl(F\cup(G+z)\bigr)
=\Es(F)+\Es(G).
\]
Taking infima at fixed volumes $m_1$ and $m_2$ yields
$e_s(m_1+m_2)\le e_s(m_1)+e_s(m_2)$.
Iterating this inequality, we find for every $N$
\[
e_s\!\left(\sum_i m_i\right)
\le \sum_{i=1}^N e_s(m_i)+e_s\!\left(\sum_{i>N}m_i\right).
\]
Since $0\le e_s(v)\le\Es(B_v)\to0$ as $v\downarrow0$
by \eqref{eq:scaling}, letting $N\to\infty$ proves \eqref{eq:countable-subadd}.
\end{proof}

\begin{lemma}\label{lem:profiles}
Let $0<s<3$, $m>0$ and $(E_j)$ be a minimizing sequence for $e_s(m)$
with $|E_j|=m$ for every $j$.
There exist a subsequence, a finite or countable family of
profiles $F^i$ of positive volume, and translations
$z_j^i\in\R^3$ such that along this subsequence
\[
 |z_j^i-z_j^{i'}|\to\infty\quad(i\ne i'),\qquad
 E_j-z_j^i\to F^i\quad\text{in }L^1_{\rm loc}(\R^3)
 \quad\text{for each }i.
\]
With $m_i=|F^i|$ one has
\begin{equation}\label{eq:profile-equality}
 \sum_i m_i=m,\qquad e_s(m)=\sum_i \Es(F^i)=\sum_i e_s(m_i).
\end{equation}
Each $F^i$ attains $e_s(m_i)$. If the family is finite, then
$|E_j\triangle\bigcup_i(F^i+z_j^i)|\to0$.
\end{lemma}

\begin{proof}
\textit{Step 1. A uniform tail estimate for the cube masses.}
For each $j$, choose distinct
$z_{j,\ell}\in\mathbb Z^3$ so that the masses
$v_{j,\ell}=|E_j\cap([0,1)^3+z_{j,\ell})|$ satisfy
\[
v_{j,1}\ge v_{j,2}\ge\cdots,
\qquad \sum_{\ell\ge1}v_{j,\ell}=m.
\]
We order the cubes of positive mass in nonincreasing order of their masses and, if there are only finitely many, complete the sequence with distinct cubes of zero mass.
By the relative isoperimetric inequality on the unit cube,
\[
\begin{aligned}
	\sum_{\ell\ge1}v_{j,\ell}^{p}
	\le& C\sum_{\ell:\,v_{j,\ell}\le1/2}
	P\bigl(E_j;(0,1)^3+z_{j,\ell}\bigr)
	+\sum_{\ell:\,v_{j,\ell}>1/2}
	v_{j,\ell}^{-1/3}v_{j,\ell}\\
	\le& CP(E_j)+2^{1/3}m,
\end{aligned}
\]
where $C$ depends only on the dimension.

After discarding finitely many terms, $P(E_j)\le\Es(E_j)\le e_s(m)+1$.
Since $(N+1)v_{j,N+1}\le\sum_{\ell=1}^{N+1}v_{j,\ell}\le m$,
\begin{equation}\label{eq:profile-tail}
	\begin{aligned}
		\sup_j\sum_{\ell>N}v_{j,\ell}
		\le& \sup_j\left(v_{j,N+1}^{1/3}
		\sum_{\ell>N}v_{j,\ell}^{p}\right)\\
		\le& \left(\frac{m}{N+1}\right)^{1/3}
		\bigl[C(e_s(m)+1)+2^{1/3}m\bigr]
		\xrightarrow[N\to\infty]{}0.
	\end{aligned}
\end{equation}
\textit{Step 2. Extraction of profiles and conservation of mass.}
By the perimeter bound in Step~1 and local $BV$ compactness, we may
pass to a common subsequence along which $E_j-z_{j,\ell}$ converges
in $L^1_{\mathrm{loc}}(\R^3)$ for every $\ell$. In particular, there exists
$v_\ell:=\lim_j v_{j,\ell}$. Since the translations lie in
$\mathbb Z^3$, a diagonal choice ensures that each difference
$z_{j,\ell}-z_{j,\ell'}$ is eventually constant or diverges in norm.
By \eqref{eq:profile-tail},
\[
0\le m-\sum_{\ell=1}^N v_\ell
=\lim_{j\to\infty}\sum_{\ell>N}v_{j,\ell}
\le\sup_j\sum_{\ell>N}v_{j,\ell}
\xrightarrow[N\to\infty]{}0.
\]

Group the indices with $v_\ell>0$ by declaring $\ell$ and $\ell'$ equivalent when
$z_{j,\ell}-z_{j,\ell'}$ remains bounded. For each class
$\mathcal I_i$, set $z_j^i=z_{j,\min\mathcal I_i}$ and denote the
local limit of $E_j-z_j^i$ by $F^i$. Then
$|z_j^i-z_j^{i'}|\to\infty$ for $i\ne i'$.

For $\ell\in\mathcal I_i$, the cube $[0,1)^3+z_{j,\ell}-z_j^i$ is eventually fixed.
By local convergence, its $F^i$-mass is $\lim_{j\to\infty}v_{j,\ell}=v_\ell$.
The corresponding lattice cubes are disjoint, so
$m_i\mathrel{=}|F^i|\ge\sum_{\ell\in\mathcal I_i}v_\ell>0$.

For any finite subfamily $1\le i\le N$ and any $R>0$, the balls
$\ball_R(z_j^i)$ are pairwise disjoint for all large $j$. Hence
\[
\sum_{i=1}^N |F^i\cap\ball_R(0)|
=\lim_{j\to\infty}
\left|E_j\cap\bigcup_{i=1}^N\ball_R(z_j^i)\right|
\le m.
\]
Letting $R\to\infty$ and then exhausting the classes, we conclude that
\[
m=\sum_i\sum_{\ell\in\mathcal I_i}v_\ell
\le\sum_i m_i\le m.
\]

\textit{Step 3. Energy equality and global approximation.}
Fix $R>0$ and finitely many profiles $F^1,\ldots,F^N$.
By the local convergence in Step~2 and \eqref{eq:riesz-L1} on any compact subinterval
of $(0,3)$ containing $s$,
$D_s(E_j\cap\ball_R(z_j^i))\to D_s(F^i\cap\ball_R(0))$
for each $i\le N$. Since the balls $\ball_R(z_j^i)$ are disjoint
for large $j$, positivity of the kernel and lower semicontinuity of
perimeter imply
\begin{equation}\label{eq:profile-energy-bound}
	\begin{aligned}
		e_s(m)
		=& \lim_{j\to\infty}\Es(E_j)\\
		\ge& \liminf_{j\to\infty}\sum_{i=1}^N
		\bigl[P(E_j;\ball_R(z_j^i))+D_s(E_j\cap\ball_R(z_j^i))\bigr]\\
		\ge& \sum_{i=1}^N
		\bigl[P(F^i;\ball_R(0))+D_s(F^i\cap\ball_R(0))\bigr].
	\end{aligned}
\end{equation}
Let $R\to\infty$ in \eqref{eq:profile-energy-bound} and take the
supremum over finite subfamilies. By \eqref{eq:countable-subadd}
and the mass identity from Step~2,
\[
e_s(m)\ge\sum_i\Es(F^i)\ge\sum_i e_s(m_i)
\ge e_s\!\left(\sum_i m_i\right)=e_s(m).
\]
This proves \eqref{eq:profile-equality}. The nonnegative differences
$\Es(F^i)-e_s(m_i)$ sum to zero, so each $F^i$ attains $e_s(m_i)$.

Suppose now that the family is finite. Its translated union has
volume at most $\sum_i m_i=m$. For fixed $R>0$, the sets
$(F^i+z_j^i)\cap\ball_R(z_j^i)$ are disjoint for large $j$.
Hence
\begin{align*}
	\limsup_{j\to\infty}
	\left|E_j\triangle\bigcup_i(F^i+z_j^i)\right|
	\le& 2m-2\sum_i\lim_{j\to\infty}
	\bigl|(E_j-z_j^i)\cap F^i\cap\ball_R(0)\bigr|\\
	=& 2\sum_i|F^i\setminus\ball_R(0)|
	\xrightarrow[R\to\infty]{}0.
\end{align*}
\end{proof}

\begin{samepage}
\subsection{Classification near the critical volume}

Comparing a pair of profiles with two balls of equal mass restricts the number of profiles at and above the critical volume.

\begingroup
\allowdisplaybreaks[0]
\interlinepenalty=10000
\begin{theorem}\label{thm:sequences}
Let $|s-1|<\varepsilon_0$, with $\varepsilon_0$ as in Theorem~\ref{thm:main},
and let $0<m\le3\mstar/2$.
Let $(E_j)$ be a minimizing sequence for $e_s(m)$
with $|E_j|=m$ for every $j$.
Along a subsequence there exist $k\in\{1,2\}$ and centers
$z_j^i\in\R^3$ such that
\[
 \left|E_j\triangle\bigcup_{i=1}^k(B_{m/k}+z_j^i)\right|
 \longrightarrow 0,\qquad |z_j^1-z_j^2|\longrightarrow\infty
 \quad\text{if }k=2.
\]
\begin{enumerate}
 \item If $m<\mstar$, then $k=1$.
 \item If $m=\mstar$, each $k\in\{1,2\}$ is realized by a minimizing sequence.
 \item If $\mstar<m\le3\mstar/2$, then $k=2$ and
       $e_s(m)=2\Es(B_{m/2})$.
\end{enumerate}
\end{theorem}
\endgroup
\end{samepage}

\begin{proof}
\textit{Step 1. Comparison of two profiles.}
By Lemma~\ref{lem:profiles} and Theorem~\ref{thm:main}, each profile
is a translate of $B_{m_i}$ with $0<m_i\le\mstar$.
Suppose there are at least two profiles. Relabel an arbitrary pair
as $F^1,F^2$ and set $v=m_1+m_2$. By \eqref{eq:countable-subadd}
and \eqref{eq:profile-equality},
\begin{equation}\label{eq:profile-pair}
	e_s(v)\le e_s(m_1)+e_s(m_2)
	=e_s(m)-\sum_{i\ge3}e_s(m_i)\le e_s(v).
\end{equation}
The strict binding inequality for $v<\mstar$
\cite[Theorem~4]{FN} therefore forces $v\ge\mstar$.

For $v\ge\mstar$, subtract \eqref{eq:gap-any-m} at $t=1/2$ from its value at $t$.
Since $G_s(1/2)=0$ and $x\mapsto x^{b_s}$ is convex,
\begin{equation}\label{eq:equal-split-comparison}
	\begin{aligned}
		&\frac{\Es(B_{tv})+\Es(B_{(1-t)v})-2\Es(B_{v/2})}{P(B_v)}\\
		=& G_s(t)+A(s)\bigl[(v/\mstar)^{q_s}-1\bigr]
		\bigl[t^{b_s}+(1-t)^{b_s}-2^{1-b_s}\bigr]
		\ge G_s(t).
	\end{aligned}
\end{equation}

For the chosen pair, take $t=m_1/v$ in \eqref{eq:equal-split-comparison}.
By \eqref{eq:profile-pair} and \eqref{eq:two-ball-upper},
\[
0\le G_s(t)
\le \frac{e_s(v)-2\Es(B_{v/2})}{P(B_v)}
\le 0.
\]
It follows from \eqref{eq:split-gap} that $t=1/2$, so
$m_1=m_2=v/2\ge\mstar/2$. Since the pair was arbitrary, all profile
masses are equal. By \eqref{eq:profile-equality}, their number $k$
is finite and $k\mstar/2\le m$.

\textit{Step 2. The number of profiles.}
Step~1 implies that $k=1$ for $m<\mstar$ and $k\in\{1,2\}$ for
$m=\mstar$. If $\mstar<m\le3\mstar/2$, then $m_i\le\mstar<m$
for every $i$, so $2\le k\le 2m/\mstar\le3$.

If $k=3$, then $m=3\mstar/2$ and all three profiles have mass $m/3$.
Since $2m/3=\mstar$, \eqref{eq:profile-equality} and the equal-split balance in
\eqref{eq:split-gap-identity} imply
\[
e_s(m)=3\Es(B_{m/3})=\Es(B_{m/3})+\Es(B_{2m/3}).
\]
Applying \eqref{eq:equal-split-comparison} at $(v,t)=(m,1/3)$
and using \eqref{eq:two-ball-upper} and \eqref{eq:split-gap}, we obtain
\[
0\ge e_s(m)-2\Es(B_{m/2})\ge P(B_m)G_s(1/3)>0.
\]
Thus $k=2$ for $\mstar<m\le3\mstar/2$, and
$e_s(m)=2\Es(B_{m/2})$ by \eqref{eq:profile-equality}.
Absorb the center of each ball $F^i$ into $z_j^i$.
The convergence and separation then follow from Lemma~\ref{lem:profiles}.

At $m=\mstar$, the constant sequence $E_j=B_m$ realizes $k=1$ by
Theorem~\ref{thm:main}. To realize $k=2$, take
$E_j=B_{m/2}\cup(B_{m/2}+z_j)$ with $|z_j|>2(R_m+R_{m/2})$ and
$|z_j|\to\infty$. By \eqref{eq:union-interaction} and
\eqref{eq:translated-interaction},
\[
\lim_{j\to\infty}\Es(E_j)=2\Es(B_{m/2})=\Es(B_m)=e_s(m).
\]
Every translate of $B_m$ intersects at most one component of $E_j$, and $B_m$ contains the component at the origin. Hence
\[
\inf_{y\in\R^3}|E_j\triangle(B_m+y)|
=2m-2\sup_{y\in\R^3}|E_j\cap(B_m+y)|=m.
\]
\end{proof}

\subsection{The exponent and separation limits}

We now let the exponent tend to $1$ and the mass to $m_c$. Along a subsequence of near-minimizers, the first-order change in the interaction is the sum of the contributions from the limiting balls.

\begin{corollary}\label{cor:varying}
Let $E_j\subset\R^3$ be measurable with
$|E_j|=m_j>0$, and let $s_j\in(0,3)$.
Assume that $s_j\to1$, $m_j\to m_c$ and
$\mathcal E_{s_j}(E_j)-e_{s_j}(m_j)\to0$. There exist a
subsequence and a fixed $k\in\{1,2\}$ such that
\[
 |E_j\triangle U_j|\to0,\qquad
 U_j=\bigcup_{i=1}^k(B_{m_j/k}+z_j^i),\qquad z_j^i\in\R^3,
\]
with $|z_j^1-z_j^2|\to\infty$ if $k=2$. Along this subsequence,
\begin{equation}\label{eq:varying}
 D_{s_j}(E_j)-D_1(E_j)
 =k\bigl(D_{s_j}(B_{m_j/k})-D_1(B_{m_j/k})\bigr)
   +|s_j-1|\,o(1).
\end{equation}
For each $M>0$, $0<v\le M$ and $s\in I$, the ball interaction
satisfies
\begin{equation}\label{eq:profile-taylor}
 \begin{aligned}
 &\biggl|D_s(B_v)-D_1(B_v)
       \mathbin{-}(s-1)D_1(B_v)
       \left(-\log 2+\frac{31}{30}-\frac13\log\frac v\omega\right)\biggr|
 \le\frac{v}{4}K_{2,M}|s-1|^2.
 \end{aligned}
\end{equation}
In particular, along the subsequence in \eqref{eq:varying},
\begin{align*}
 D_{s_j}(E_j)-D_1(E_j)
 =& (s_j-1)kD_1(B_{m_c/k})
       \left(-\log 2+\frac{31}{30}-\frac13\log\frac{m_c}{k\omega}\right)
       +|s_j-1|\,o(1).
\end{align*}
\end{corollary}

\begin{proof}
\textit{Step 1. Reduction to a Coulomb minimizing sequence.}
Fix $M>m_c$. For all large $j$, we have $s_j\in I$ and $m_j\le M$.
Set $\tau_j=(m_c/m_j)^{1/3}$, so that $\tau_j\to1$ and
$|\tau_jE_j|=m_c$. By \eqref{eq:scaling}, \eqref{eq:threshold},
and Theorem~\ref{thm:main},
\[
\mathcal E_{s_j}(B_{m_j})
\longrightarrow \mathcal E_1(B_{m_c})=e_1(m_c).
\]
Applying \eqref{eq:exponent-lipschitz} to $\tau_jE_j$ and using
near minimality, we obtain
\begin{equation}\label{eq:coulomb-reduction}
	\begin{aligned}
		e_1(m_c)
		\le& \mathcal E_1(\tau_jE_j)\\
		\le& \tau_j^2P(E_j)+\tau_j^{6-s_j}D_{s_j}(E_j)
		+\frac{m_c}{2}K_{1,M}|s_j-1|\\
		\le& \max\{\tau_j^2,\tau_j^{6-s_j}\}\,
		\mathcal E_{s_j}(B_{m_j})+o(1)
		\longrightarrow e_1(m_c).
	\end{aligned}
\end{equation}

Thus $(\tau_jE_j)$ is a minimizing sequence for $e_1(m_c)$ by \eqref{eq:coulomb-reduction}.
By Theorem~\ref{thm:sequences}, a subsequence has $k\in\{1,2\}$ ball profiles of mass $m_c/k$. Write the centers
of the approximating balls as $\tau_jz_j^i$. Since
$\tau_jB_{m_j/k}=B_{m_c/k}$,
\[
|E_j\triangle U_j|
=\frac{m_j}{m_c}
\left|\tau_jE_j\triangle
\bigcup_{i=1}^k(B_{m_c/k}+\tau_jz_j^i)\right|
\longrightarrow0.
\]
If $k=2$, the rescaled centers separate and $\tau_j\to1$, so
$|z_j^1-z_j^2|\to\infty$. Their distance eventually exceeds
$\max\{4R_{M/2},2\}$, and the two balls are then disjoint.
Thus $|U_j|=m_j$ for all large $j$.

\textit{Step 2. First-order change in the interaction.}
Since $|E_j\triangle U_j|\to0$, \eqref{eq:parameter-L1} implies
\begin{equation}\label{eq:varying-decomposition}
	\bigl|(D_{s_j}-D_1)(E_j)-(D_{s_j}-D_1)(U_j)\bigr|
	\le K_{1,M}|s_j-1|\,|E_j\triangle U_j|
	=|s_j-1|\,o(1).
\end{equation}
For $k=1$, $U_j$ is a translate of $B_{m_j}$, so \eqref{eq:varying} follows immediately.
For $k=2$, we estimate the change in the interaction between the two balls.
For $z\in\R^3$ with $|z|\ge2$, let $F,G\subset\ball_{|z|/4}(0)$ be measurable.
Then
$1\le |z|/2\le |x-y|\le 2|z|$ whenever $x\in F$ and $y\in G+z$.
For $s\in I$, the mean value theorem in the exponent yields
\begin{equation}\label{eq:cross-small}
	\begin{aligned}
		&\left|\int_F\int_{G+z}
		\bigl(|x-y|^{-s}-|x-y|^{-1}\bigr)\dd y\dd x\right|\\
		\le& |s-1|\int_F\int_{G+z}
		|x-y|^{-s_-}\log|x-y|\dd y\dd x\\
		\le& 2^{s_-}|F|\,|G|\,|s-1|\,
		\frac{\log(2|z|)}{|z|^{s_-}}.
	\end{aligned}
\end{equation}
Taking $F=G=B_{m_j/2}$ and $z=z_j^2-z_j^1$ in \eqref{eq:cross-small} and using \eqref{eq:union-interaction}, we obtain
\begin{align*}
	&\bigl|(D_{s_j}-D_1)(U_j)-2(D_{s_j}-D_1)(B_{m_j/2})\bigr|\\
	\le&
	|s_j-1|\,
	\frac{m_j^2\log\!\bigl(2|z_j^1-z_j^2|\bigr)}
	{2^{2-s_-}|z_j^1-z_j^2|^{s_-}}
	=|s_j-1|\,o(1),
\end{align*}
since $m_j\to m_c$, $|z_j^1-z_j^2|\to\infty$ and $r^{-s_-}\log(2r)\le 2^{1+s_-/2}r^{-s_-/2}/s_-$ for $r\ge1$.
Thus \eqref{eq:varying} follows from \eqref{eq:varying-decomposition}.

For $|E|\le M$ and $s\in I$, Taylor's formula
at $s=1$ and \eqref{eq:kernel-derivatives} imply
\begin{equation}\label{eq:uniform-taylor}
  \left|D_s(E)-D_1(E)
    +\frac{s-1}{2}\int_E\int_E\frac{\log|x-y|}{|x-y|}\dd y\dd x\right|
  \le \frac{|E|}{4}K_{2,M}|s-1|^2 .
\end{equation}
For $v>0$, differentiating the ball formula in \eqref{eq:scaling}
and using \eqref{eq:unit-interaction-log-derivative}, we find
\[
\left.\partial_s D_s(B_v)\right|_{s=1}
=D_1(B_v)\left(\frac{31}{30}-\log 2
-\frac13\log\frac{v}{\omega}\right).
\]
Taking $E=B_v$ in \eqref{eq:uniform-taylor} proves
\eqref{eq:profile-taylor} for $0<v\le M$.

Since the balls are concentric,
$k|B_{m_j/k}\triangle B_{m_c/k}|=|m_j-m_c|$.
Applying \eqref{eq:parameter-L1} and \eqref{eq:profile-taylor}
to the ball term in \eqref{eq:varying}, we conclude that
\begin{align*}
	&\left|(D_{s_j}-D_1)(E_j)
	-k(s_j-1)\bigl[\partial_s D_s(B_{m_c/k})\bigr]_{s=1}\right|\\
	\le&|s_j-1|
	\left(K_{1,M}|m_j-m_c|
	+\frac{m_c}{4}K_{2,M}|s_j-1|+o(1)\right)
	=|s_j-1|\,o(1).
\end{align*}
\end{proof}

\begin{samepage}
	\section*{Statements and Declarations}
	
	\paragraph{Funding}
	The research of Chengxiang Zhang was supported by the
	National Natural Science Foundation of China
	(Grant No. 12371107).
	
	\paragraph{Data availability}
	No datasets were generated or analysed in this study.
	
	\paragraph{Use of generative AI}
	The authors acknowledge the use of ChatGPT (OpenAI) to assist with
	selected calculations. They have independently verified all
	mathematical arguments and calculations and take full responsibility
	for the content of the paper.
	
	
\end{samepage}

\begingroup
\renewcommand{\bibfont}{\small\raggedright}
\setlength{\bibsep}{3pt plus 1pt minus 1pt}

\endgroup
\end{document}